\documentclass[12pt, dvips, twoside]{amsart}
\usepackage{amssymb,latexsym,bbm,graphicx,epsfig,epic,eepic,oldgerm}
\usepackage{eucal, amscd}
\usepackage{psfrag}
\usepackage{a4wide}
                                    
\theoremstyle{plain}
\newtheorem{theorem}{Theorem}[section]
\newtheorem{lemma}[theorem]{Lemma}
\newtheorem{proposition}[theorem]{Proposition}
\newtheorem{corollary}[theorem]{Corollary}
\newtheorem{remark}[theorem]{Remark}
\newtheorem{remarks}[theorem]{Remarks}

\newtheorem{examples}[theorem]{Examples}

\newtheorem{question}[theorem]{Question}
\newtheorem{NStheorem}[theorem]{Nonsqueezing Theorem}

\newcommand{\proofend}{\hspace*{\fill} $\Box$\\}

\def\s{\smallskip}
\def\m{\medskip}

\def\eps{\varepsilon}

\def\cat{\operatorname{cat}}

\def\cl{\operatorname{cl}}
\def\B{\operatorname {B}}
\def\C{\operatorname {C}}
\def\SB{\operatorname {S_B}}
\def\SBe{\operatorname {S_B^=}}
\def\SBa{\mbox{S}_{\mbox{\tiny B}}^a} 
\def\S{\operatorname {S}}
\def\Sdis{\operatorname {S_{dis}}}
\def\Gr{\operatorname {Gr\;\!}}
\def\dim{\operatorname {dim}}

\def\Vol{\operatorname {Vol}\:\!}
\def\Int{\operatorname {Int}\:\!}
\def\dist{\operatorname {dist}\:\!}
\def\inter{\operatorname {int}}
\def\ext{\operatorname {ext}}
\def\diameter{\operatorname {diam}\:\!}
\def\Emb{\operatorname {Emb}}
\def\id{\operatorname {id}}

\def\gg{\gamma}
\def\gd{\delta}
\def\eps{\varepsilon}

\def\gf{\varphi}

\def\gl{\lambda}
\def\go{\omega}
\def\gs{\sigma}
\def\gt{\vartheta}

\def\ca{{\mathcal A}}
\def\cb{{\mathcal B}}
\def\cc{{\mathcal C}}

\def\cg{{\mathcal G}}

\def\ck{{\mathcal K}}
\def\cn{{\mathcal N}}

\def\cp{{\mathcal P}}

\def\cs{{\mathcal S}}
\def\ct{{\mathcal T}}
\def\cu{{\mathcal U}}

\def\cw{{\mathcal W}}

\def\ea{{\mathfrak A}}

\def\ec{{\mathfrak C}}
\def\ed{{\mathfrak D}}

\def\eq{{\mathfrak Q}}
\def\es{{\mathfrak S}}

\def\CC{\mathbbm{C}}

\def\NN{\mathbbm{N}}
\def\PP{\mathbbm{P}}

\def\RR{\mathbbm{R}}

\def\ZZ{\mathbbm{Z}}
\def\CP{{\CC\PP}} 

\def\11{\mathbbm{1}}

\def\pp{\partial}

\def\ra{\rightarrow}
\def\ha{\hookrightarrow}
\def\Ra{\Rightarrow}

\def\ni{\noindent}
\def\b{\bigskip}
\def\m{\medskip}

\def\proof{\noindent {\it Proof. \;}}

\begin{document}



\title[]{Minimal atlases of closed symplectic manifolds}

\date{\today}

\author{Yu.~B.~Rudyak}
\address{(Yu.~B.~Rudyak) University of Florida, 
Department of Math., 358 Little Hall,
PO Box 118105, Gainesville, FL 32611-8105, USA}
\email{rudyak@math.ufl.edu}
\author{Felix Schlenk}
\address{(F.~Schlenk)         
D\'epartement de Math\'ematiques, Universit\'e Libre de Bruxelles, 
CP~218, Boulevard du Triomphe, 1050 Bruxelles, Belgium}
\email{fschlenk@ulb.ac.be}

\begin{abstract}  
We study the number of Darboux charts needed to cover a closed
connected symplectic manifold $(M,\go)$ and effectively estimate 
this number from below and from above in terms of
the Lusternik--Schnirelmann category of $M$ and the Gromov width of $(M,\go)$.
\end{abstract}

\maketitle

\section{Introduction and main results}  \label{s:intro}

\ni
A symplectic manifold is a pair $(M, \go)$ where $M$ is a
smooth manifold and $\go$ is a non-degenerate and closed $2$-form on $M$. 
The non-degeneracy of $\go$ implies that $M$ is even-dimensional,
$\dim M = 2n$.
(We refer to \cite{HZ} and \cite{MS} for basic facts about symplectic
manifolds.)
The most important symplectic manifold is $\RR^{2n}$ equipped with its
standard symplectic form
\[
\go_0 \,=\, \sum_{i=1}^n dx_i \wedge dy_i .
\]
Indeed,
a basic fact about symplectic manifolds is Darboux's Theorem which
states that locally every symplectic manifold $(M^{2n}, \go)$ is diffeomorphic 
to $(\RR^{2n}, \go_0)$.
More precisely, for each point $p \in M$ there exists a chart 
\[
\gf \colon B^{2n} (a) \ra M
\]
from a ball
\[
B^{2n} (a) :\,= \left\{ z \in \RR^{2n} \mid \pi \left| z \right|^2 < a
\right\}
\]
to $M$ such that $\gf (0) = p$ and $\gf^* \go = \go_0$.
We call such a chart $\left( B^{2n}(a), \gf \right)$ a Darboux chart. 
In this paper we study the following question:

\m
\begin{quote}
{\it Given a closed symplectic manifold $(M, \go)$, how many Darboux
charts does one need in order to parametrize $(M, \go)$?}
\end{quote}
\m

\ni
In other words, we study the number $\SB (M ,\go)$ defined as 
\[
\SB (M, \go) \,:\,=\, \min \left\{ k \mid M = \cb_1 \cup \dots \cup \cb_k 
\right\}
\]
where each $\cb_i$ is the image $\gf_i \left( B^{2n}(a_i) \right)$ of a
Darboux chart.
 
\m
An obvious lower bound for $\SB (M, \go)$ is the diffeomorphism
invariant 
\[
\B (M) \,:\,=\, \min \left\{ k \mid M =B_1 \cup \dots \cup B_k \right\}
\]
where each $B_i$ is diffeomorphic to the standard open ball in $\RR^{2n}$. 

The volume associated with a symplectic manifold $\left( M^{2n}, \go
\right)$ is
\[
\Vol (M, \go) \,=\, \frac{1}{n!} \int_M \go^n .
\]
In particular, $\Vol \left( B^{2n}(a) \right) = \tfrac{1}{n!} \;\!
a^n$, as it should be.
The volume of any symplectically embedded ball in $(M, \go)$ is at most 
\[
\gg \left( M, \go\right)  \,=\, 
\sup \left\{ \Vol \left( B^{2n}(a) \right) \mid B^{2n}(a) 
\text{ symplectically embeds into } M \right\} .
\]
Another lower bound for $\SB (M, \go)$ is therefore
\[
\Gamma (M, \go) \,:\,=\, \left\lfloor \frac{\Vol (M, \go)}{\gg (M, \go)}
\right\rfloor +1 
\]
where $\lfloor x \rfloor$ denotes the maximal integer which is
smaller than or equal to $x$.
Notice that $\gg (M, \go) = \tfrac{1}{n!} \left( \Gr (M, \go) \right)^n$ 
where
\[
\Gr (M,\go) = \sup \left\{ a \mid B^{2n}(a) \text{ symplectically embeds
into } (M, \go) \right\} 
\]
is the Gromov width of $(M, \go)$. The symplectic invariant $\Gamma (M,
\go)$ is therefore strongly related to the Gromov width. 
We abbreviate
\[
\gl (M, \go) \,:\,=\, \max \left\{ \B(M),\:\! \Gamma (M, \go) \right\} .
\]
Summarizing we have that
\begin{equation}  \label{e:lsb}
\gl (M, \go) \,\le\, \SB (M, \go) .
\end{equation}
Before we state our main result, we consider two examples.

\s
1)
For complex projective space $\CP^n$ equipped with its standard
K\"ahler form
$\go_{SF}$ we have $\B (\CP^n) = n+1$ and $\Gamma (\CP^n, \go_{SF}) = 
2$.
In particular,
\[
\gl (\CP^n, \go_{SF}) \,=\, \B (\CP^n) \,>\, \Gamma (\CP^n,
\go_{SF}) \quad \text{if } n \ge 2 .
\]
It will turn out that 
$\SB (\CP^n, \go_{SF}) = \gl (\CP^n, \go_{SF}) = n+1$ if $n \ge 2$. 

\s
2)
We fix an area form $\gs$ on the 2-sphere $S^2$, and for $k \in \NN$
we abbreviate $S^2(k) = (S^2, k \gs)$.
Then
$\B \left( S^2 \times S^2 \right) = 3$ and $\Gamma \left( S^2(1) \times
S^2(k) \right) = 2k+1$.
In particular,
\[
\gl \left( S^2(1) \times S^2(k) \right)  \,=\, 
\Gamma \left( S^2(1) \times S^2(k) \right) \,>\, 
\B \left( S^2 \times S^2 \right) \quad \text{if } k \ge 2 .
\]
It will turn out that $\SB \left( S^2(1) \times S^2(k) \right) = 
\gl \left( S^2(1) \times S^2(k) \right) = 2k+1$ if $k \ge 2$. 

\s
We refer to Examples 2 and 4 in Section \ref{s:ex} for more details.

\m
\ni
Our main result is

\m
\ni
{\bf Theorem~1.}
{\it Let $(M, \go)$ be a closed connected $2n$-dimensional 
symplectic manifold.
\begin{itemize}
\item[(i)]
If $\gl (M, \go) \ge 2n+1$, then $\SB (M, \go) = \gl (M, \go)$.
\s
\item[(ii)]
If $\gl (M, \go) < 2n+1$, then $n+1 \le \gl (M, \go) \le \SB (M, \go)
\le 2n+1$. 
\end{itemize}
}

\b
\ni
{\bf Remarks.}
{\bf 1.}
The assumption in (i) is met if $[\go] |_{\pi_2(M)} =0$, see
Proposition~1\:(ii) below.
It is also met for various symplectic fibrations, see Section~\ref{s:ex}.

\s
{\bf 2.}
Theorem 1 implies that
\[
n+1 \,\le\, \gl (M, \go) \,<\, \SB (M,\go) \,\le\, 2n+1 
\quad \text{ if }\;
\gl (M, \go) \neq \SB (M,\go) .
\]

\s
\ni
The following question is based on the examples described in Section~\ref{s:ex}.

\m
\ni
{\bf Question.}
{\it
Is it true that $\gl (M, \go) = \SB (M,\go)$ for all closed symplectic
manifolds $(M, \go)$?
}

\b
Theorem~1 essentially reduces the problem of computing the number
$\SB(M,\go)$ to two other problems, namely computing $\B (M)$ and
$\Gamma (M, \go)$. 
The computation of the Gromov width and hence of $\Gamma (M, \go)$ 
is often a very delicate matter. 
Fortunately, there has recently been remarkable
progress in this problem, 
see \cite{B1,B2,BC,Ji,KT,LM1,LM3,Lu,M-variants,MP,MSl,Sch1,Sib} and
Section~\ref{s:ex}.
On the other hand, the diffeomorphism invariant $\B(M)$ can often 
be computed or estimated very well,
as we shall explain next.

Recall that the {\it Lusternik--Schnirelmann category} of a finite
$CW$-space $X$ is defined as
\[
\cat X \,:\,=\, \min \left\{ k \mid X = A_1 \cup \ldots \cup A_k \right\} , 
\]
where each $A_i$ is open and contractible in $X$, \cite{LS, CLOT}.
Clearly,
\[
\cat M \,\le\, \B (M) 
\]
if $M$ is a closed smooth manifold.
It holds that
$\cat X =\cat Y$ whenever $X$ and $Y$ are homotopy equivalent. However, the  
Lusternik--Schnirelmann
category is very different from the usual homotopical invariants
in algebraic topology and hence often difficult to compute.
Nevertheless, $\cat X$ can be estimated from below in cohomological terms as
follows. 
Let $H^*$ be singular cohomology, with any coefficient ring, 
and let $\tilde{H}^*$ be the corresponding reduced cohomology.
The {\it cup-length} of $X$ is 
defined as
\[
\cl(X) \,:\,=\, \sup \left\{ k \mid u_1 \cdots u_k \neq 0, u_i \in
\tilde{H}^*(X) \right\}.
\]
It then holds true that
\[
\cl (X) +1 \le \cat X,
\]
see \cite{FE}. 
Much more information on LS-category can be found in \cite{{CLOT}, {J1}, {J2}}.

If $M^m$ is a smooth closed connected manifold, then 
$\B (M) \le m+1$, see~\cite{L,Z}. 
Summarizing we have that 
\begin{equation} \label{e:F}
\cl (M)+1 \le \cat M \le \B(M) \le m+1 
\end{equation}
for any closed $m$-dimensional manifold. 

\medskip
These inequalities may be
substantially improved if $M$ is symplectic.    

\m
\ni
{\bf Proposition~1.}  
{\it Let $(M, \go)$ be a closed connected $2n$-dimensional 
symplectic manifold.
Then
\begin{equation}  \label{p:prop1}
n+1 \le \cl (M) +1 \le \cat M \le \B(M) \le 2n+1.
\end{equation}
Moreover, the following assertions hold true.
\begin{itemize}
\item[(i)]
If $\pi_1(M) =0$, then $n+1 = \cl(M)+1 = \cat M = \B(M)$.
\s
\item[(ii)]
If $[\omega] |_{\pi_2(M)} =0$, then $\cat M = \B(M) = 2n+1$.
\s
\item[(iii)]
If $\cat M < \B(M)$, then $n \ge 2$, $n+1 = \cl(M) +1 = \cat M$ and
$\B(M) = n+2$. 

\end{itemize}
}

\b
The paper is organized as follows.
In Section~2 we prove Theorem~1.
In Section~3 we study the minimal number $\SBe (M,\go)$ of {\it equal}\,
symplectic balls needed to cover $(M,\go)$
as well as the minimal number $\S (M,\go)$ of symplectic charts diffeomorphic
to a ball needed to parametrize $(M,\go)$.
In Section~4 we prove Proposition~1, 
and in the last section we compute the number $\SB (M,\go)$
for various closed symplectic manifolds.

\b
\ni
{\bf Outlook.}
In a sequel, we shall use the symplectic ball covering number $\SB$ to formulate 
a Lusternik--Schnirelmann theory for (Wein-)Stein manifolds and polarized 
K\"ahler manifolds as studied in 
\cite{E-A,EG} and \cite{B3,BC}.

\b
\ni
{\bf Acknowledgments.}
The idea behind the proof of Theorem~1 belongs to Gromov. We are very
grateful to Leonid Polterovich for explaining it to us. The first author was 
supported by NSF, grant 0406311.

\section{Proof of Theorem~1}  \label{s:t1}

\ni
In view of the inequalities \eqref{e:lsb} and \eqref{p:prop1}, Theorem~1 is a
consequence of

\begin{theorem}  \label{t:ref}
Let $(M, \go)$ be a closed connected $2n$-dimensional symplectic manifold.
\begin{itemize}
\item[(i)]
If $\Gamma (M, \go) \ge 2n+2$, then $\SB (M,\go) = \Gamma (M, \go)$.

\s
\item[(ii)]
If $\Gamma (M, \go) \le 2n+1$, then $\SB (M,\go) \le 2n+1$.
\end{itemize}
\end{theorem} 

\b
\ni
{\bf Idea of the proof.}
We start with describing the idea of the proof, which belongs to
Gromov and is as simple as beautiful.
For each Borel set $A$ in $M$ we abbreviate its volume
\[
\mu (A) \,:\,=\, \frac{1}{n!} \int_A \omega^n .
\]
Moreover, we define the natural number $k$ by
\begin{equation}  \label{def:k}
k \,=\, 
\left\{
 \begin{array}{lll} 
    \Gamma (M, \go) &\text{ if} & \Gamma (M, \go) \ge 2n+2, \\
    2n+1            &\text{ if} & \Gamma (M, \go) \le 2n+1.
 \end{array}
\right.
\end{equation} 
By definition of $\Gamma (M, \go)$,
\begin{equation}  \label{est:gmk}
\gg (M, \go) \,>\, \frac{\mu (M)}{k} .
\end{equation}
By definition of $\gg (M, \go)$ we find a Darboux chart $\gf \colon
B^{2n}(a) \ra \cb \subset M$ such that
\[
\mu \left( \cb \right) \,>\, \frac{ \mu (M)}{k} .
\]
In view of this inequality, and since $k \ge 2n+1 = \dim M+1$, 
elementary dimension theory will provide
a cover of $M$ by $k$ sets $\cc^1, \dots, \cc^k$ where each set $\cc^j$ is
essentially a disjoint union of small cubes, and where
\[
\mu \left( \cc^j \right) \,<\, \mu \left( \cb \right)
\quad \text{for each } j,
\]
cf.~Figure~\ref{figure20} below.
Using this and the specific choice of the sets $\cc^j$ we shall then be
able to construct for each $j$ a symplectomorphism $\Phi^j$ of $M$
such that $\Phi^j(\cc^j) \subset \cb$.
The $k$ Darboux charts
\[
\left( \Phi^j \right)^{-1} \circ \gf \colon B^{2n}(a) \ra M
\]
will then cover $M$, and so Theorem~\ref{t:ref} follows.

\begin{figure}[ht] 
 \begin{center}
  \psfrag{1}{$\Phi^j$}
  \leavevmode\epsfbox{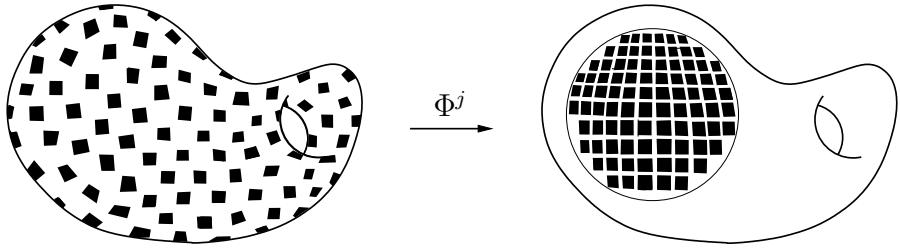}
 \end{center}
 \caption{The idea behind the map $\Phi^j$.} 
 \label{figure24}
\end{figure}
%
%

\ni
Notice that $\mu \left( \cc^j \right)$ might be very close to $\mu
\left( \cb \right)$.
In order that the ``cubes'' in $\cc^j$ all fit into the ball $\cb$, 
the map $\Phi^j$ should therefore not distort the cubes too much. 
We shall be able to find such a map $\Phi^j$ by constructing an
appropriate atlas for $(M,\go)$ and by constructing the set $\cc^j$ carefully.

\b
\ni
{\bf Step~1. Construction of a good atlas of $(M, \go)$} 

\begin{figure}[ht] 
 \begin{center}
  \psfrag{p1}{$p_1$}
  \psfrag{p2}{$p_2$}
  \psfrag{U0}{$\cu_0$}
  \psfrag{U0'}{$\cu_0'$}
  \psfrag{U1}{$\cu_1$}
  \psfrag{U1'}{$\cu_1'$}
  \psfrag{U2}{$\cu_2$}
  \psfrag{U2'}{$\cu_2'$}
  \leavevmode\epsfbox{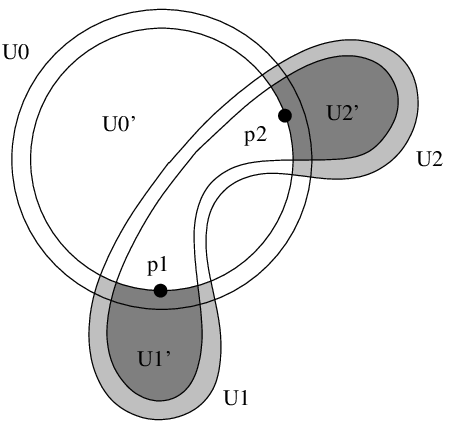}
 \end{center}
 \caption{The sets $\cu_1' \subset \cu_1$ and $\cu_2' \subset \cu_2$
   and the points $p_1 \in \pp \cu_0' \cap \pp \cu_1'$ and 
                  $p_2 \in \pp \cu_0' \cap \pp \cu_2'$.}  
 \label{figure27}
\end{figure}

\s
\ni
Let $k$ be the natural number defined in \eqref{def:k}. In view of
the estimate \eqref{est:gmk} the real number $\eps$ defined by 
\[
\gg (M,\go) \,=\, \frac{\mu(M)}{k} + 2 \eps
\]
is positive.
By definition of $\gg (M, \go)$ we can choose 
a Darboux chart 
\[
\gf_0 \colon B^{2n}(a_0) \ra \cb_0 \subset M
\]
such that 
\begin{equation*}
\mu (\cb_0) \,>\, \frac{\mu(M)}{k} + \eps .
\end{equation*}
Since $M$ is compact, we find $m$ other Darboux charts
$\gf_i \colon B^{2n}(a_i) \ra \cb_i \subset M$ such that
\begin{equation}  \label{id:MB}
M \,=\, \bigcup_{i=0}^m \cb_i .
\end{equation}
We can assume that 
\begin{equation}  \label{inc:nt}
\cb_i \,\not\subset\, \bigcup_{j \neq i} \cb_j, \quad \, i=0,\dots,m .
\end{equation}
Given open subsets $U \subset V$ of $\RR^{2n}$ we write $U \Subset V$
if $\overline{U} \subset V$,
and we say that a symplectic chart $\bigl( \widetilde{U}, \widetilde{\gf} 
\bigr)$ is {\it larger}\, 
than a symplectic chart $\left( U, \gf \right)$ if $U \Subset
  \widetilde{U}$ and $\gf = \widetilde{\gf} |_U$. 
Using this terminology we can also assume that each chart 
$\left( B^{2n}(a_i), \gf_i \right)$ is the restriction of a larger chart. 
Then the boundaries of the images $\cb_0, \cb_1, \dots, \cb_m$ are smooth.
%
%
We next choose for $i=0, \dots, m$ numbers $a_i' < a_i$ so large that
with $\cb_i' = \gf_i \left( B^{2n}(a_i') \right)$ we have
\begin{equation}  \label{est:B0p}
\mu \left( \cb_0' \right) \,>\, \frac{\mu(M)}{k} + \eps 
\end{equation}
and
\begin{equation}  \label{id:prime}
M \,=\, \bigcup_{i=0}^m \cb_i' .
\end{equation}
After renumbering the charts 
$\left( B^{2n}(a_1), \gf_1 \right),\dots, \left( B^{2n}(a_m), \gf_m \right)$
we can then assume that $\cb_1 \cap \cb_0' \neq \emptyset$. 
In view of \eqref{inc:nt} and since the boundaries of $\cb_1$ and $\cb_0'$
are smooth, the open set
\[
\cb_1 \setminus \overline{\cb_0'} \,=:\, \coprod_{i=1}^{I_1} \cu_i
\]
is non-empty and consists of finitely many connected components $\cu_i$ with
piecewise smooth boundaries.
For notational convenience we set $\cu_0 = \cb_0$ and $\cu_0' = \cb_0'$ as
well as 
\[
\cu_i' \,=\, \cu_i \cap \cb_1', \quad \, i=1, \dots, I_1.
\]
After choosing $a_1' < a_1$ larger if necessary, 
we can assume that each $\cu_i'$ is non-empty and also connected.
Each $\cu_i'$ has piecewise smooth boundary $\pp \cu_i'$.
Clearly,
\begin{equation}  \label{id:u2}
\bigcup_{i=0}^1 \cb_i \,=\, \bigcup_{i=0}^{I_1} \cu_i 
\qquad \text{and} \qquad
\bigcup_{i=0}^1 \overline{\cb_i'} \,=\, \bigcup_{i=0}^{I_1} \overline{\cu_i'},
\end{equation}
cf.~Figure~\ref{figure27}.

%
%

\ni
For each $i \in \{ 1, \dots, I_1 \}$ we choose a point
\[
p_i \,\in\, \pp \cb_0' \cap \pp \cu_i' .
\]
We let $\ct_1$ be the rooted tree whose vertices are the root $p_0$ and
the points $p_i$ and whose edges are $\left[ p_0, p_i \right]$, 
$i = 1, \dots, I_1$.
The tree $\ct_1$ corresponding to Figure~\ref{figure27} is depicted in
Figure~\ref{figure23}.
We also set $U_0 = B^{2n}(a_0)$ and $\phi_0 = \gf_0 \colon U_0 \ra
\cu_0$ and define the symplectic charts 
\[
U_i = \gf_1^{-1} \left( \cu_i \right), \quad
\phi_i = \gf_1 |_{U_i} \colon U_i \ra \cu_i , \quad \, i =1, \dots, I_1.
\]
Notice that each chart $(U_i, \phi_i)$ is the restriction of a
larger chart.

\begin{figure}[ht] 
 \begin{center}
  \psfrag{U}{$\cu_3$}
  \psfrag{U'}{$\cu_3'$}
  \psfrag{p3}{$p_3$}
  \leavevmode\epsfbox{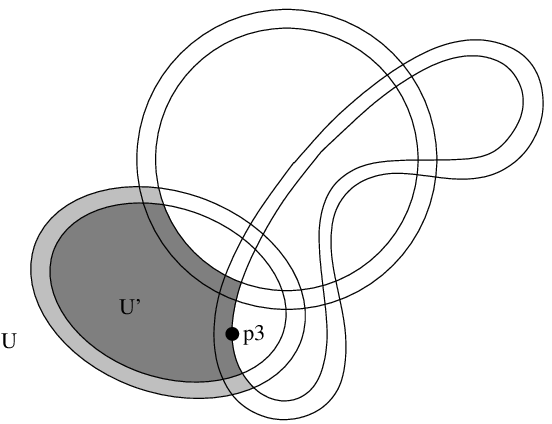}
 \end{center}
 \caption{The sets $\cu_3' \subset \cu_3$
   and the point $p_3 \in \pp \cu_1' \cap \pp \cu_3'$.}  
 \label{figure29}
\end{figure}

\begin{figure}[ht] 
 \begin{center}
  \psfrag{01}{$p_0$}
  \psfrag{11}{$p_1$}
  \psfrag{12}{$p_2$}
  \psfrag{21}{$p_3$}
  \leavevmode\epsfbox{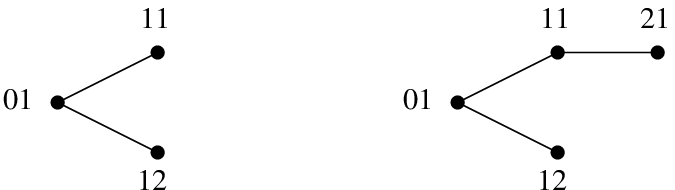}
 \end{center}
 \caption{The trees $\ct_1$ and $\ct_2$.} 
 \label{figure23}
\end{figure}

If $m \ge 2$, assumption~\eqref{id:prime} implies that we can 
renumber the charts $\left( B^{2n}(a_2), \gf_2 \right),
\dots, \left( B^{2n}(a_m), \gf_m \right)$
such that $\cb_2 \cap \bigcup_{i=0}^1 \cb_i' \neq \emptyset$.
In view of \eqref{inc:nt} and since the boundaries of $\cb_2$, $\cb_0'$
and $\cb_1'$ are smooth, the open set
\begin{equation}  \label{def:Vi}
\cb_2 \setminus \bigcup_{i=0}^1 \overline{\cb_i'} \,=:\, 
\coprod_{i=I_1+1}^{I_2} \cu_i 
\end{equation}
is non-empty and consists of finitely many connected components $\cu_i$ with
piecewise smooth boundaries.
We set $\cu_i' = \cu_i \cap \cb_2'$ for $i = I_1+1, \dots, I_2$.
After choosing $a_2' < a_2$ larger if necessary, 
each $\cu_i'$, $i = I_1+1, \dots, I_2$, is non-empty and connected,
and has piecewise smooth boundary.
Clearly,
\[
\bigcup_{i=0}^2 \cb_i \,=\, \bigcup_{i=0}^{I_2} \cu_i 
\qquad \text{and} \qquad
\bigcup_{i=0}^2 \overline{\cb_i'} \,=\, \bigcup_{i=0}^{I_2}
\overline{\cu_i'} ,
\]
cf.~Figure~\ref{figure29}.

%
%

\ni
In view of the second identity in~\eqref{id:u2} and the 
definition~\eqref{def:Vi} of $\cu_i$ we find for each 
$i \in \left\{ I_1+1, \dots, I_2 \right\}$ an index 
$\underline{i} \in \left\{ 0, \dots, I_1 \right\}$ 
such that $\pp \cu_{\underline{i}}' \cap \pp \cu_i' \neq \emptyset$,
and we choose a point
\[
p_i \,\in\, \pp \cu_{\underline{i}}' \cap \pp \cu_i' .
\]
We let $\ct_2$ be the tree obtained from the tree $\ct_1$ by adding
the vertices $p_i$ and the edges 
$\left[ p_{\underline{i}}, p_i \right]$, $i = I_1+1, \dots, I_2$. 
The tree $\ct_2$ corresponding to Figure~\ref{figure29} is depicted in 
Figure~\ref{figure23}.

%
%

\ni
We define the symplectic charts
\[
U_i = \gf_2^{-1} \left( \cu_i \right), \quad
\phi_i = \gf_2 |_{U_i} \colon U_i \ra \cu_i , \quad \, i =I_1+1, \dots, I_2.
\]
Notice again that each chart $(U_i, \phi_i)$ is the restriction of a larger 
chart.

Proceeding this way $m-2$ other times we find a sequence
\[ 
0 =: I_0 < I_1 < \dots < I_m =: l
\]
of integers and $l+1$ open connected sets $\cu_i \subset M$, 
$i=0, \dots, l$, with piecewise smooth boundaries such that for each
$j \in \{ 0, \dots, m-1 \}$, 
\begin{equation}  \label{id:BBV}
\cb_{j+1} \setminus \bigcup_{i=0}^j \overline{\cb_i'} \,=:\, 
\coprod_{i=I_j+1}^{I_{j+1}} \cu_i .
\end{equation}
Moreover, defining $j(i)$ by the condition 
$i \in \left\{ I_{j(i)}+1, \dots, I_{j(i)+1} \right\}$,  
we see that each set $\cu_i' = \cu_i \cap \cb_{j(i)+1}'$ is non-empty and 
connected, and has piecewise smooth boundary.
Furthermore, we have found for each 
$i \in \left\{ 1, \dots, l \right\}$ an index 
$\underline{i} \in \{ 0, \dots, I_{j(i)} \}$ 
such that $\pp \cu_{\underline{i}}' \cap \pp \cu_i' \neq \emptyset$
and have chosen a point
\begin{equation}  \label{inc:pvv}
p_i \,\in\, \pp \cu_{\underline{i}}' \cap \pp \cu_i' .
\end{equation}
The vertices of the rooted tree $\ct = \ct_m$ consist of the root $p_0$
and the points $p_i$, and the edges of $\ct$ are 
$\left[ p_{\underline{i}}, p_i \right]$, 
$i =1, \dots, l$. 

In view of~\eqref{id:BBV},
\begin{equation}  \label{eq:uiuj}
\cu_i' \cap \cu_j = \emptyset \quad \text{ if } \; i<j.
\end{equation}
Moreover, the identities~\eqref{id:MB} and \eqref{id:BBV} imply that
\begin{equation}  \label{id:MV}
M \,=\, \bigcup_{i=0}^l \cu_i 
\end{equation}
and that
$\sum_{i=0}^l \mu \left( \cu_i \right) \ra \mu (M)$ as $a_j'
\ra a_j$ for all $j =0,\dots,m$.
Choosing $a_0', \dots, a_m'$ larger if necessary we can therefore
assume that
\begin{equation}  \label{est:mme}
\sum_{i=0}^l \mu \left( \cu_i \right) \,<\, \mu (M) + \eps .
\end{equation}

We replace the symplectic atlas 
$\left\{ \gf_i \colon B^{2n}(a_i) \ra \cb_i, \: i=0,\dots,m  \right\}$
by the symplectic atlas
$\left\{ \phi_i \colon U_i \ra \cu_i, \: i=0,\dots,l \right\}$.
Here, we still have 
$\left( U_0, \phi_0 \right) = \left( B^{2n}(a_0), \gf_0 \right)$, and 
\[
U_i = \gf_{j(i)+1}^{-1} \left( \cu_i \right), \quad
\phi_i = \gf_{j(i)+1} |_{U_i} \colon U_i \ra \cu_i ,
\quad \,
i = 1, \dots, l .
\]
Each chart $(U_i, \phi_i)$ is the restriction of a larger chart 
$\widetilde{\phi}_i \colon \widetilde{U}_i \ra \widetilde{\cu}_i$. 
While $p_i \notin \cu_i$ in view of \eqref{inc:pvv},
we have $p_i \in \cu_{\underline{i}} \cap
\widetilde{\cu}_i$ for $i = 1, \dots, l$.

Our next goal is to replace the charts $\widetilde{\phi}_i \colon
\widetilde{U}_i \ra \widetilde{\cu}_i$ by charts
$\widetilde{\psi}_i \colon \widetilde{V}_i \ra \widetilde{\cu}_i$ 
such that for each $i \ge 1$ the transition function
\[
\widetilde{\psi}_{\underline{i}}^{-1} \circ \widetilde{\psi}_i \colon \;
\widetilde{\psi}_i^{-1} \bigl(
  \widetilde{\cu}_{\underline{i}} \cap \widetilde{\cu}_i \bigr)
\,\ra\,
 \widetilde{\psi}_{\underline{i}}^{-1} \bigl( \widetilde{\cu}_{\underline{i}} 
\cap
 \widetilde{\cu}_i \bigr) 
\]
is the identity on a neighbourhood $W_i$ of 
$\widetilde \psi_i^{-1} (p_i)$.
The neighbourhoods $\cw_i = \widetilde \psi_i (W_i)$ will serve as gates
for moving cubes from $\widetilde \cu_i$ to $\widetilde \cu_{\underline i}$ 
without distorting them.
We first of all set 
$\bigl( \widetilde{V}_0, \widetilde{\psi}_0 \bigr) =
\bigl( \widetilde{U}_0, \widetilde{\phi}_0 \bigr)$.
In order to construct $\bigl( \widetilde{V}_1, \widetilde{\psi}_1 \bigr)$
we first define a symplectic chart 
$\bigl( \widehat{V}_1, \widehat{\psi}_1 \bigr)$ by
\[
\widehat{V}_1 = \left[ d \left( \widetilde{\phi}_1^{-1} \circ
    \widetilde{\psi}_0 \right) (q_1) \right]^{-1} \bigl(
  \widetilde{U}_1 \bigr), 
\quad
\widehat{\psi}_1 = 
\widetilde{\phi}_1 \circ d \left( \widetilde{\phi}_1^{-1} \circ
    \widetilde{\psi}_0 \right) (q_1) \colon \widehat{V}_1 \ra
  \widetilde{U}_1 
\]
where we abbreviated $q_1 = \widetilde{\psi}_0^{-1} (p_1)$.
We then find
\begin{equation}  \label{id2:psi}
\left( \widetilde{\psi}_0^{-1} \circ \widehat{\psi}_1\right) (q_1) = q_1
\quad
\text{ and }
\quad
d \left( \widetilde{\psi}_0^{-1} \circ \widehat{\psi}_1 \right) (q_1)
= \id . 
\end{equation}
We obtain the desired chart 
$\bigl( \widetilde{V}_1, \widetilde{\psi}_1 \bigr)$
from the chart  
$\bigl( \widehat{V}_1, \widehat{\psi}_1 \bigr)$ with the help of the
following lemma.

\begin{lemma}  \label{l:normal}
Assume that $\gf \colon U \ra U'$ is a symplectomorphism between two
domains $U$ and $U'$ in $\RR^{2n}$ such that $\gf (q) =q$ and 
$d\gf(q) = \id$ at some point $q \in U$.
Then there exist open neighbourhoods $W \subset \widetilde{W} \Subset
U$ of $q$ 
and a symplectomorphism $\rho \colon U \ra U'$ such that 
$\rho |_W = id$ and 
$\rho |_{U \setminus \widetilde{W}} = \gf |_{U \setminus
  \widetilde{W}}$.  
\end{lemma}

\proof
We can assume that $q=0$.
Following \cite[Appendix A.1]{HZ} we represent the map $\gf$ by
\begin{eqnarray*}
x &=& a(\xi, \eta) \\
y &=& b(\xi, \eta) .
\end{eqnarray*}
Since $d \gf (0) = \id$, we have $\det \left( a_{\xi} (0) \right) =1
\neq 0$.
According to Proposition~1 in \cite[Appendix A.1]{HZ} we
therefore find a smooth function $w$ defined on a neighbourhood $\cn
\subset \RR^{2n}(x, \eta)$ of $0$ such that 
\begin{eqnarray}  \label{e:xyw}
\left\{
 \begin{array}{lcl} 
    \xi &=& x+ w_{\eta} (x,\eta) \\
    y   &=& \eta + w_x (x,\eta) .
 \end{array}
\right.
\end{eqnarray}
We can assume that $w(0) =0$.
In view of the identities $\gf(0) =0$ and $d \gf (0) = \id$ and the
relations \eqref{e:xyw} we find that all the derivatives of $w$ up to
order $2$ vanish in $0$, i.e., 
\begin{equation}  \label{id:wo}
w(x,\eta) \,=\, O \left( \left| (x,\eta) \right|^3 \right) .
\end{equation}
Choose a smooth function $f \colon [0, \infty[ \;\ra [0,1]$ such that
\[ 
f(s) \,=\,
\left\{
 \begin{array}{ll} 
    0, & s \le 1, \\
    1, & s \ge 2,
 \end{array}
\right.
\]
and denote the open ball of radius $s$ in $\RR^{2n}(x,\eta)$ by $B_s$.
For each $\eps >0$ for which $B_{3\eps} \subset \cn$ we define the
smooth function $w^{\eps} (x,\eta) \colon B_{3\eps} \ra \RR$ by
\[
w^{\eps} (x,\eta) \,=\, f \left( \tfrac{1}{\eps} \left| (x,\eta) \right|
\right) w(x,\eta) .
\]
Then
\begin{eqnarray}  \label{id:weu}
w^{\eps} |_{B_{\eps}} =0 
\qquad \text{and} \qquad
w^{\eps} |_{B_{3\eps}\setminus B_{2\eps}} = w |_{B_{3\eps}\setminus
  B_{2\eps}} . 
\end{eqnarray}
Abbreviating $\zeta = (x,\eta)$ and $r = \left| \zeta \right|$ 
we compute
\begin{eqnarray*}
w^{\eps}_{\zeta_i}(\zeta) &=& 
  f' \left( \tfrac{r}{\eps} \right) \tfrac{1}{\eps} \:\!
  \tfrac{\zeta_i}{r}  \;\!
  w(\zeta) + f \left( \tfrac{r}{\eps} \right) w_{\zeta_i}(\zeta), \\
w^{\eps}_{\zeta_i \zeta_j}(\zeta) &=&
  f'' \left( \tfrac{r}{\eps} \right) \tfrac{1}{\eps^2} \:\! \tfrac{\zeta_i
    \zeta_j}{r^2} \;\! w(\zeta) 
  +  f' \left( \tfrac{r}{\eps} \right) \tfrac{1}{\eps} \left(
    \tfrac{\gd_{ij}}{r} - \tfrac{\zeta_i \zeta_j}{r^3} \right) w(\zeta)
                                                                     \\
 & &  \qquad \qquad \qquad \quad \;\;\:
        + \:f' \left( \tfrac{r}{\eps} \right) \tfrac{1}{\eps}
 \left( \tfrac{\zeta_i}{r} \;\! w_{\zeta_j} (\zeta) + \tfrac{\zeta_j}{r}
                                   \;\! w_{\zeta_i} (\zeta)  \right)
                                      \\  
 & &  \qquad \qquad \qquad \quad \;\;\:
        + \:f \left( \tfrac{r}{\eps} \right) w_{\zeta_i \zeta_j}(\zeta) 
\end{eqnarray*}
where $i,j \in \{1, \dots, 2n\}$ and where $\gd_{ij}$ denotes the Kronecker
symbol. In view of the estimate~\eqref{id:wo} we therefore find that 
\[
w^{\eps}_{\zeta_i \zeta_j}(\zeta) \,=\, \tfrac{1}{\eps^2} O (r^3) +
              \tfrac{1}{\eps} O (r^2) + O(r) \,=\, 
              O(r), \quad \zeta \in B_{3\eps},   
\]
and so 
\begin{eqnarray}  \label{est:wij}
w^{\eps}(x,\eta) \,=\, O \left( \left| (x,\eta) \right|^3 \right), 
                            \quad \, (x,\eta) \in B_{3\eps}.        
\end{eqnarray}
We in particular conclude that 
$\det \bigl( \11_n + w^\eps_{x\eta}(x,\eta) \bigr) 
\neq 0$ for all $(x,\eta) \in B_{3\eps}$ if $\eps >0$ is small enough.
The relations
\begin{eqnarray}  \label{e:xiy}
\left\{
 \begin{array}{lcl} 
    \xi &=& x+ w_{\eta}^{\eps} (x,\eta) \\
    y   &=& \eta + w_x^{\eps} (x,\eta) 
 \end{array}
\right.
\end{eqnarray}
therefore implicitly define a symplectic mapping 
$\gf^{\eps} \colon (\xi, \eta) \mapsto (x,y)$ near $0$, see again \cite[Appendix 
A.1]{HZ}.
The $C^2$-estimate \eqref{est:wij} 
implies that $\gf^{\eps}$ is $C^1$-close to the identity and that for
$\eps >0$ small enough, $\gf^{\eps}$ is defined and injective on all
of 
\[
U_{3\eps}^{\eps} \,=\, \left\{ (\xi, \eta) \in \RR^{2n} \mid
  \text{\eqref{e:xiy} holds for } (x,\eta) \in B_{3\eps} \right\} .
\]
In view of the estimate \eqref{est:wij} each of the sets 
\[
U_s^{\eps} \,=\, \left\{ (\xi, \eta) \in \RR^{2n} \mid
  \text{\eqref{e:xiy} holds for } (x,\eta) \in B_s \right\}, \quad
s \le 3 \eps,  
\]
is contained in the domain $U$ of $\gf$ and is diffeomorphic to an open ball 
provided
that $\eps >0$ is small enough.
According to the identities~\eqref{id:weu}, 
the map $\gf^{\eps}$ is the identity on $U_{\eps}^{\eps}$ and
coincides with $\gf$ on the ``open annulus'' $U_{3\eps}^{\eps}
\setminus \overline{U_{2\eps}^{\eps}}$.
It follows that $\gf^{\eps} \left( U_{3\eps}^{\eps} \right) = \gf
\left( U_{3\eps}^{\eps} \right)$.
We smoothly extend $\gf^{\eps} \colon U_{3\eps}^{\eps} \ra \RR^{2n}$
to a symplectic embedding $\rho \colon U \ra \RR^{2n}$ by setting
$\rho (z) = \gf(z)$, $z \in U \setminus U_{3\eps}^{\eps}$.
Then $\rho (U) = \gf (U) = U'$, and setting 
$W = U_{\eps}^{\eps}$ and 
$\widetilde{W} = U_{2\eps}^{\eps} \Subset U_{3\eps}^{\eps} \subset U$ 
we find that 
$\rho |_W = \gf^{\eps} |_{U_{\eps}^{\eps}} = id$ and 
$\rho |_{U \setminus \widetilde{W}} = \gf |_{U \setminus \widetilde{W}}$.
The proof of Lemma~\ref{l:normal} is complete.
\proofend

In view of the identities~\eqref{id2:psi} we can apply Lemma~\ref{l:normal} to 
the
symplectomorphism 
\[
\widetilde{\psi}_0^{-1} \circ \widehat{\psi}_1 \colon \;
\widehat{\psi}_1^{-1} \bigl( \widetilde{\cu}_0 \cap \widetilde{\cu}_1
\bigr) \,\ra\,  
\widetilde{\psi}_0^{-1} \bigl( \widetilde{\cu}_0 \cap \widetilde{\cu}_1
\bigr)
\]
which fixes $q_1$, and find open neighbourhoods $W_1 \subset
\widetilde{W}_1 \Subset \widehat{\psi}_1^{-1} \bigl( \widetilde{\cu}_0
  \cap \widetilde{\cu}_1 \bigr)$ and a symplectomorphism
\[
\rho_1 \colon \,
\widehat{\psi}_1^{-1} \bigl( \widetilde{\cu}_0 \cap \widetilde{\cu}_1 \bigr) 
\,\ra\,
\widetilde{\psi}_0^{-1} \bigl( \widetilde{\cu}_0 \cap \widetilde{\cu}_1 \bigr) 
\]
such that
\begin{equation}  \label{prop:r}
\rho_1 |_{W_1} = id 
\quad \text{ and } \quad
\rho_1 |_{\widehat{\psi}_1^{-1} \left( \widetilde{\cu}_0 \cap
    \widetilde{\cu}_1 \right) \setminus \widetilde{W}_1} 
= \widetilde{\psi}_0^{-1} \circ \widehat{\psi}_1 .
\end{equation}
Set $\widetilde{V}_1 = \widehat{V}_1$.
In view of the properties~\eqref{prop:r} of $\rho_1$ the map
$\widetilde{\psi}_1 \colon \widetilde{V}_1 \ra \widetilde{\cu}_1$
defined by
\[
\widetilde{\psi}_1 \,=\,                       
    \left\{
    \begin{array}{lll}   
        \widetilde{\psi}_0 \circ \rho_1 & \text{on} &
          \widehat{\psi}_1^{-1} \bigl( \widetilde{\cu}_0 \cap
          \widetilde{\cu}_1 \bigr),   \\ [0.3em]
        \widehat{\psi}_1 & \text{on} &\widetilde{V}_1 \setminus
        \widetilde{W}_1  
    \end{array}   \right. 
\] 
is a well-defined smooth symplectic chart such that
\[
\widetilde{\psi}_0^{-1} \circ \widetilde{\psi}_1 \colon \;
\widetilde{\psi}_1^{-1} \bigl( \widetilde{\cu}_0 \cap \widetilde{\cu}_1
\bigr) \,\ra\,  
\widetilde{\psi}_0^{-1} \bigl( \widetilde{\cu}_0 \cap \widetilde{\cu}_1
\bigr)
\]
is the identity on the open neighbourhood $W_1$ of
$q_1=\widetilde{\psi}_0^{-1} (p_1)$.
Assume now by induction that we have already constructed new charts 
$\widetilde{\psi}_j \colon \widetilde{V}_j \ra \widetilde{\cu}_j$ for
$j = 1, \dots, i-1$.
Since $\underline{i} < i$, the chart 
$\bigl( \widetilde{U}_{\underline{i}},
\widetilde{\phi}_{\underline{i}} \bigr)$ 
is already replaced by the chart
$\bigl( \widetilde{V}_{\underline{i}}, \widetilde{\psi}_{\underline{i}}
\bigr)$.
Applying the two-step construction exemplified above to the pair   
$\bigl( \widetilde{V}_{\underline{i}}, \widetilde{\psi}_{\underline{i}}
\bigr)$,
$\bigl( \widetilde{U}_{i}, \widetilde{\phi}_{i}
\bigr)$  
we find a new chart
$\widetilde{\psi}_{i} \colon \widetilde{V}_{i} \ra
\widetilde{\cu}_{i}$ such that the transition function
\[
\widetilde{\psi}_{\underline{i}}^{-1} \circ \widetilde{\psi}_{i} \colon \;
\widetilde{\psi}_i^{-1} \bigl(
  \widetilde{\cu}_{\underline{i}} \cap \widetilde{\cu}_{i} \bigr)
\,\ra\,
 \widetilde{\psi}_{\underline{i}}^{-1} \bigl(
   \widetilde{\cu}_{\underline{i}} \cap \widetilde{\cu}_{i} \bigr) 
\]
is the identity on an open neighbourhood $W_{i}$ of 
$q_{i} = \widetilde{\psi}_{\underline{i}}^{-1} (p_{i})$.
In this way we obtain a new symplectic atlas
\[
\widetilde{\ea} \,=\, \left\{ \widetilde{\psi}_i \colon \widetilde{V}_i 
\,\ra\,
\widetilde{\cu}_i, \: i = 0, \dots, l \right\} .
\]
Recall that $\cu_i \Subset \widetilde{\cu}_i$.
The collection 
\[
\ea \,=\, \left\{ \psi_i \colon V_i 
\,\ra\,
\cu_i, \: i = 0, \dots, l \right\} 
\]
of smaller charts defined by
\[
V_i = \widetilde{\psi}_i^{-1} \left( \cu_i \right), \quad
\psi_i = \widetilde{\psi}_i |_{V_i} \colon V_i \ra \cu_i
\]
is the good atlas of $(M, \go)$ we were looking for. 
For later reference we summarize the properties of this atlas:
\begin{itemize}
\item[1.]
The chart $\psi_0 \colon V_0 \ra \cu_0$ is equal to $\gf_0 \colon B^{2n}(a_0) 
\ra \cb_0$.

\s
\item[2.]
For each $i=0, \dots, l$ the chart $\psi_i \colon V_i \ra \cu_i$
is the restriction of a larger chart $\widetilde \psi_i \colon \widetilde V_i 
\ra \widetilde \cu_i$.
Each set $\cu_i$ is connected and has piecewise smooth boundary, and 
contains a certain domain $\cu_i'$ with piecewise smooth boundary.

\s
\item[3.]
There is a rooted tree $\ct$ whose root corresponds to $\cu_0$, 
whose vertices correspond to $\cu_0, \dots,\cu_l $,
and whose edges correspond to points 
$p_i \in \pp \cu_{\underline i}' \cap \pp \cu_i'$ where $i=1, \dots, l$ and 
$\underline i < i$.
Each $p_i$ has an open neighbourhood $\cw_i \subset \cu_{\underline i}' \cap 
\cu_i'$ on which the transition function 
$\widetilde \psi_{\underline i}^{-1} \circ \widetilde \psi_i$ is the identity.
\end{itemize}

\b
\ni
{\bf Step~2. The dimension cover $\ed (2n,k)$} 

\s
\ni
Let $k \ge 2n+1$ be the natural number defined in \eqref{def:k}.
In this step we construct a special cover $\ed (2n,k)$ of
$\RR^{2n}$ by cubes.
Our construction is inspired by an idea from elementary dimension
theory, see e.g.~\cite[Figure~7]{E}.

We denote the coordinates in $\RR^{2n}$ by $x_1,\dots, x_{2n}$,
and we let $\left\{ e_1, \dots, e_{2n} \right\}$ be the standard basis
of $\RR^{2n}$. 
Given a point $q \in \RR^{2n}$ and a subset $A$ of $\RR^{2n}$ we
denote the translate of $A$ by $q$ by
\[
q+A \,=\, \left\{\:\! q+a \mid a \in A \:\!\right\} .
\]
By a cube we mean a translate of the closed cube 
$C^{2n} = [0,1]^{2n} \subset \RR^{2n}$.
We define the $\left( 2n \times 2n \right)$-matrix $M(2n,k)$ as the
matrix whose diagonal is $(k,1, \dots,1)$, whose upper-diagonal is 
\[
\left( \tfrac{k}{2n}, \tfrac{2n}{2n-1}, \tfrac{2n-1}{2n-2}, \dots,
  \tfrac{4}{3}, \tfrac{3}{2} \right)
\]
and whose other matrix entries are zeroes.
E.g., 
\[
M(2,3) = \begin{bmatrix} 3 & \frac{3}{2} \\
                         0 & 1 \end{bmatrix},
\quad
M(2,4) = \begin{bmatrix} 4 & 2 \\
                         0 & 1 \end{bmatrix},
\quad
M(4,5) = \begin{bmatrix} 5 & \frac{5}{4} &0  &  0 \\
                         0 & 1 & \frac{4}{3} & 0 \\
                         0 & 0 & 1 & \frac{3}{2} \\
                         0 & 0 & 0 & 1   \end{bmatrix}.
\]
We consider the infinite union of cubes
\[
\ec^1 (2n,k) \,=\, \bigcup_{v \in \ZZ^{2n}} M(2n,k) v + C^{2n}
\]
and its translates
\[
\ec^j (2n,k) \,=\, (j-1) e_1 + \ec^1(2n,k), \quad j = 2, \dots, k,
\]
and we define the cover
\[
\ed (2n,k) \,:\,=\, \left\{ \ec^j (2n,k) \right\}_{j=1}^k,
\]
cf.~Figure~\ref{figure20} and Figure~\ref{figure21}.

\begin{figure}[ht] 
 \begin{center}
  \psfrag{1}{$x_1$}
  \psfrag{2}{$x_2$}
  \leavevmode\epsfbox{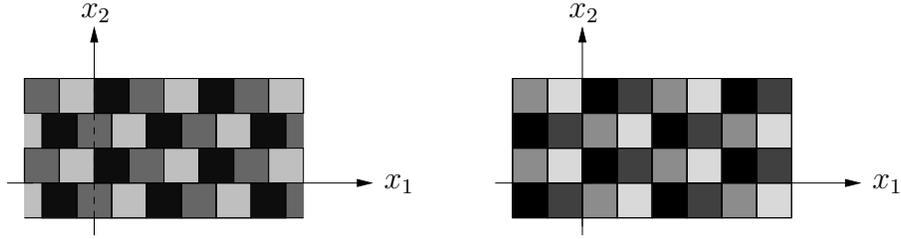}
 \end{center}
 \caption{Parts of the dimension covers $\ed (2,3)$ and $\ed (2,4)$.} 
 \label{figure20}
\end{figure}
%
%

%
%
%
\begin{figure}[ht] 
 \begin{center}
  \psfrag{x1}{$x_1$}
  \psfrag{x2}{$x_2$}
  \leavevmode\epsfbox{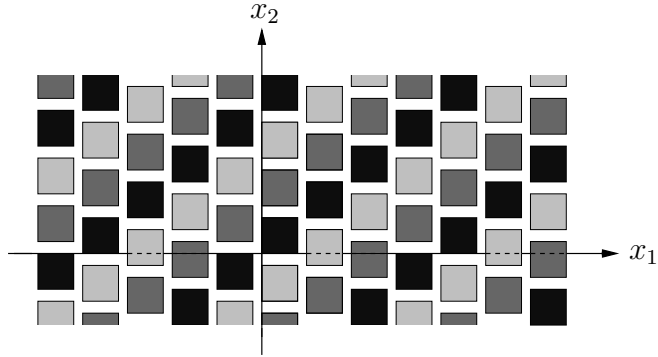}
  \end{center}
\caption{A part of the intersections
   $\ec^1 (4,5) \cap 
   \{ (x_1, x_2,
   x_3, x_4) \mid x_3 = i - \tfrac{1}{2}, x_4 =0 \}$, $i=1,2,3$.} 
\label{figure21}
\end{figure}
%
%

\ni
Finally, we define for each subset $A$ of $\RR^{2n}$ and each $m \in \{ 1,
\dots, 2n \}$ the cylinder $Z_m(A)$ over $A$ by
\[
Z_m(A) \,=\, \left\{ a + \gl e_m \mid a \in A, \, \gl \in \RR \right\} .
\]
Recall that the distance between two subsets $A$ and $B$ of
$\RR^{2n}$ is defined as
\[
\dist(A,B) \,=\, \inf \left\{ \left| a-b \right| \mid a \in A, \, b \in B
 \right\} .
\]
Given $\nu >0$ and a subset $A$ of $\RR^{2n}$ we denote the
$\nu$-neighbourhood of $A$ by
\[
\cn_{\nu}(A) \,=\, \left\{ z \in \RR^{2n} \mid \dist (z,A) < \nu 
                                                        \right\} .
\]
We abbreviate the positive number
\begin{equation}  \label{def:dd}
\gd \,:\,=\, \min \left\{ \tfrac{k-2n}{2n}, \tfrac{1}{2n-1} \right\} .
\end{equation}

\begin{lemma}\  \label{l:cover}

\begin{itemize}
\item[(i)]
For each $j \in \{1,\dots,k\}$ and any cube $C$ of $\ec^j(2n,k)$ we have
\[
\dist (C,\ec^j(2n,k) \setminus C) \,=\, \gd .
\]
Moreover,
\[
Z_1 \left( \Int C \right) \cap \ec^j(2n,k) \,=\, \bigcup_{l\in\ZZ} kle_1 + \Int 
C
\]
and 
\[
Z_m \left( \cn_{\gd} (C) \right) \cap \ec^j(2n,k) \,=\, \bigcup_{l\in\ZZ} 
(2n-m+2)le_m + C, \quad
m =2,\dots,2n.
\]
\item[(ii)] The family $\ed (2n,k)$ is a cover of $\RR^{2n}$, i.e., 
\[
\bigcup_{j=1}^k \ec^j(2n,k) \,=\, \RR^{2n} ,
\]
and the interiors of the sets $\ec^j(2n,k)$ are mutually disjoint.
\end{itemize}
\end{lemma} 
\ni
The proof, which is elementary, is omitted.

\b
\ni
{\bf Step~3. The cover of $M$ by small cubes} 

\s
\ni
Let $\ea = \left\{ \psi_i \colon V_i \ra \cu_i, \: i=0, \dots, l
                                                         \right\}$  
be the symplectic atlas of $(M, \go)$ constructed in Step~1 and
let $\ed (2n, k) = \left\{ \ec^j(2n,k)\right\}_{j=1}^k$ be the dimension
cover of $\RR^{2n}$ constructed in the previous step.
For any $r>0$ and any subset $A$ of $\RR^{2n}$ we set
\[
rA \,=\, \left\{ rz \mid z \in A \right\} 
\]
and we denote by $\left| A \right|$ the Lebesgue measure of $A$.
Fix $i \in \{ 0, \dots, l \}$.
For $d_i >0$ we define $\ec_i^j (d_i)$ as the union of those cubes $C$ in 
$d_i \ec^j (2n,k)$ for which  
\begin{equation}  \label{cond}
C \subset V_i \quad \text{ and } \quad
\dist \left( C, \pp V_i \right) \ge d_i
\end{equation}
and we abbreviate
\[
\ed_i (d_i) \,:\,=\, \bigcup_{j=1}^k \ec_i^j (d_i) .
\]
By ``a cube of $\ec_i^j(d_i)$" we shall mean a component of $\ec_i^j(d_i)$,
and 
by ``a cube of $\ed_i(d_i)$" we shall mean a cube of some $\ec_i^j(d_i)$.
In view of the identity~\eqref{id:MV} 
we find open sets $\breve{\cu}_i \Subset \cu_i$ such that 
\[
M \,=\, \bigcup_{i=0}^l \cu_i \,=\, \bigcup_{i=0}^l \breve{\cu}_i .
\]
Choose $d_i>0$ so small that $\psi_i^{-1} \bigl( \breve{\cu}_i \bigr)
\subset \ed_i(d_i)$. Then
\begin{equation}  \label{id:MU}
M \,=\, \bigcup_{i=0}^l \psi_i \bigl( \ed_i (d_i) \bigr) .
\end{equation}
Also notice that the ``homogeneity'' of the sets $\ec_i^j(d_i)$ implies
that
\[
\left| \ec_i^j(d_i) \right| \,\ra\, \tfrac{1}{k} \left| V_i \right|
\quad \text{ as }\;  d_i \ra 0
\]
for all $j \in \{ 1, \dots, k \}$. Choosing $d_i >0$ smaller if
necessary we can therefore assume that
\begin{equation}  \label{e:Cijdi}
\left| \ec_i^j(d_i) \right| \,<\, 
\tfrac{1}{k} \left( \left| V_i \right| + \tfrac{k-1}{l+1} \eps \right)
\end{equation}
for all $i \in \{ 0, \dots, l \}$ and $j \in \{ 1, \dots, k \}$.

\begin{figure}[ht] 
 \begin{center}
  \leavevmode\epsfbox{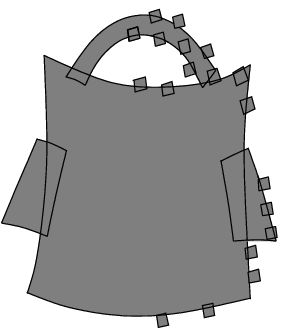}
 \end{center}
 \caption{A component of $\cc_2^j$.} 
 \label{figure.height}
\end{figure}

We denote by $\cc^j = \cc^j (d_0, \dots, d_l)$ the union of cubes ``of
the same colour $j$'', 
\[
\cc^j \,=\, \bigcup_{i=0}^l \psi_i \bigl( \ec_i^j (d_i) \bigr), \quad
\, j = 1, \dots, k.
\]
The components $\psi_i (C)$ of  $\psi_i \bigl( \ec_i^j (d_i) \bigr)$ are
called $i$-cubes.
For each connected component $\ck$ of $\cc^j$ we define the {\it
  height} of $\ck$ as the maximal $h \in \{ 0, \dots, l \}$ for which
$\ck$ contains an $h$-cube. 
The set $\cc^j$ decomposes as 
\[
\cc^j \,=\, \coprod_{h=0}^l \cc_h^j 
\]
where $\cc_h^j$ is the union of the components of $\cc^j$ of height $h$,
cf.~Figure~\ref{figure.height}.

%

\ni
In view of~\eqref{id:MU} we have
\begin{equation}  \label{id:MUU}
M \,=\, \bigcup_{j=1}^k \bigcup_{h=0}^l \cc_h^j .
\end{equation}
According to the estimates \eqref{e:Cijdi} we can choose for each $i
\in \left\{ 1, \dots, l \right\}$ a number
\begin{equation}  \label{est:nud}
\nu_i \,\in\, \left] 0, \tfrac{\gd}{2} \right[
\end{equation}
such that
\begin{equation}  \label{est:2nu}
\left(1+ 2\nu_i\right)^{2n} \left| \ec_i^j (d_i) \right|  \,<\, 
   \tfrac{1}{k} \left( \left| V_i \right| + \tfrac{k-1}{l+1} \eps \right)
\end{equation}
for all $j \in \left\{ 1, \dots, k \right\}$.
Since $\nu_i < \tfrac{\gd}{2} < 1$, the conditions \eqref{cond}
imply that
\begin{equation}  \label{inc:NV}
\cn_{\nu_i d_i} (C) \,\subset\, V_i
\end{equation}
for any cube $C$ of $\ed_i (d_i)$.

\begin{lemma}  \label{l:nu}
If the numbers $d_0, \dots, d_{l-1} >0$ as well as the ratios 
$d_i / d_{i+1}$, $i = 0, \dots, l-1$, are small enough,
then the following assertions hold true.
\begin{itemize}
\item[(i)]
$\cc_h^j \subset \cu_h$\, for each $j \in \{ 1, \dots,k \}$ and $h \in
\{ 0, \dots, l \}$.
\item[(ii)]
Any component $\ck$ of $\cc_h^j$ contains only one $h$-cube $\psi_h(C)$,
and 
\[
\psi_h^{-1} (\ck) \,\subset\, \cn_{\nu_h d_h}(C) , \quad \, h =1, \dots,
l.
\]
\end{itemize}
\end{lemma} 

\proof
We denote by $\cp_i^j = \cp_i^j (d_0, \dots, d_l)$ the partial union 
of cubes 
\[
\cp_i^j \,=\, \bigcup_{g=i}^l \psi_i \bigl( \ec_i^j (d_i) \bigr), \quad
\, i=0, \dots,l; \;\, j = 1, \dots, k.
\]
E.g., $\cp_l^j = \psi_l \bigl( \ec_l^j (d_l) \bigr)$ and $\cp_0^j = \cc^j$.
Generalizing the above definition we define the 
{\it  height} of a connected component $\ck$ of $\cp^j_i$ 
as the maximal $h \in \{ i, \dots, l \}$ for which
$\ck$ contains an $h$-cube. 
The set $\cp_i^j$ decomposes as 
\[
\cp_i^j \,=\, \coprod_{h=i}^l \cp_{i,h}^j 
\]
where $\cp_{i,h}^j$ is the union of components of $\cp_i^j$ of height $h$.

Since $\cp_l^j$ consists of finitely many disjoint closed cubes, we can choose
$d_{l-1} > 0$ so small that each cube of
$\psi_{l-1} \bigl( \ec_{l-1}^j \left(d_{l-1} \right) \bigr)$ intersects
at most one cube of $\cp^j_l$ for each $j$. 
Then each component $\ck$ of $\cp_{l-1,l}^j$ contains only one
$l$-cube.
We denote the distinguished cube in $\ck$ by $\cc (\ck)$.
Since $\cp_l^j$ is a compact subset of the open set $\cu_l$, we can
choose $d_{l-1}$ so small that $\cp_{l-1,l}^j \subset \cu_l$ for each
$j$. Moreover, choosing $d_{l-1}$ yet smaller if necessary we can
assume that
\begin{equation}  \label{inc:plk}
\psi_l^{-1} (\ck) \,\subset\, \cn_{\nu_ld_l} \left( \psi_l^{-1} \left(
    \cc (\ck) \right) \right) 
\end{equation}
for each component $\ck$ of $\cp_{l-1,l}^j$ and each $j$.

Since $\cp_{l-1}^j$ consists of finitely many disjoint compact
components, we can choose $d_{l-2}>0$ so small that each cube of
$\psi_{l-2} \bigl( \ec_{l-2}^j \left(d_{l-2} \right) \bigr)$ intersects
at most one component of $\cp^j_{l-1}$ for each $j$. 
Then each component $\ck$ of $\cp_{l-2,h}^j$ contains only one
$h$-cube, $h = l,l-1,l-2$.
We denote this distinguished cube again by $\cc (\ck)$.
If $h \in \left\{ l,l-1 \right\}$, then $\cc (\ck) = \cc \left(
  \underline{\ck} \right)$ where $\underline{\ck}$ is the unique
component of $\cp_{l-1,h}^j$ contained in $\ck$, and if $h = l-2$, then
$\cc (\ck) = \ck$ is an $(l-2)$-cube.
Since $\cp_{l-1,l}^j$ is a compact subset of the open set $\cu_l$ and
since $\cp^j_{l-1, l-1}$ is a compact subset of the open set
$\cu_{l-1}$, 
we can choose $d_{l-2}$ so small that $\cp_{l-2,l}^j \subset \cu_l$ 
and $\cp_{l-2,l-1}^j \subset \cu_{l-1}$ for each $j$.
Moreover, the compact inclusions~\eqref{inc:plk} imply that we can choose
$d_{l-2}$ so small that
\[
\psi_l^{-1} (\ck) \,\subset\, \cn_{\nu_ld_l} \left( \psi_l^{-1} \left(
    \cc (\ck) \right) \right) 
\]
for each component $\ck$ of $\cp_{l-2,l}^j$ and each $j$.
Choosing $d_{l-2}$ yet smaller if necessary we can also
assume that
\[
\psi_{l-1}^{-1} (\ck) \,\subset\, \cn_{\nu_{l-1}d_{l-1}} 
  \left( \psi_l^{-1} \left( \cc (\ck) \right) \right) 
\]
for each component $\ck$ of $\cp_{l-2,l-1}^j$ and each $j$.


Repeating this reasoning $l-2$ other times, we successively find
$d_{l-1}, \dots, d_0$ such that assertions (i) and (ii) of the lemma
hold true for all $h \in \left\{ 1, \dots, l \right\}$ and all
$j$. Since $\cc_0^j \subset \cu_0$ by definition of $\cc_0^j$, the
proof of Lemma~\ref{l:nu} is complete.
\proofend

\begin{figure}[ht] 
 \begin{center}
  \leavevmode\epsfbox{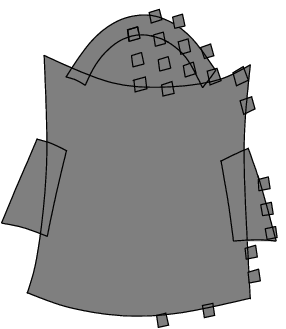}
 \end{center}
 \caption{A component of $\cs_2^j$.} 
 \label{figure28}
\end{figure}

For $h \ge 1$ the sets $M \setminus \cc_h^j$ are not necessarily 
connected.
We define the {\em saturation}\, $\cs (A)$ of a compact subset $A$
of $\RR^{2n}$ as the union of $A$ with the bounded components of
$\RR^{2n} \setminus A$.
Since $A$ is compact, $\RR^{2n} \setminus \cs (A)$ is the only unbounded 
component of $\RR^{2n} \setminus A$ and hence in particular is connected.
For a compact subset $\ca$ of $\cu_h$ with $\cs \left( \psi_h^{-1}
  \left( \ca \right) \right) \subset V_h$ we define its saturation as
\[
\cs (\ca) \,=\, \psi_h \left( \cs \left( \psi_h^{-1} (\ca) \right)
\right) .
\]
By Lemma~\ref{l:nu}\:(ii) and the inclusions~\eqref{inc:NV} we have
$\cs \bigl( \psi_h^{-1} \bigl( \cc_h^j \bigr) \bigr) \subset V_h$ for all
$j \in \left\{ 1, \dots, k \right\}$ and $h \in \left\{ 0, \dots, l
\right\}$. For $j \in \left\{ 1, \dots, k \right\}$ we can therefore
recursively define compact subsets of $\cu_h$ by
\begin{eqnarray*}
\cs_l^j  &=& \cs \left( \cc_l^j \right), \\
\cs_h^j  &=& \cs \left( \cc_h^j \setminus \bigcup_{g=h+1}^l \cs_g^j
                      \right),  \quad \; h = l-1, \dots, 0 .
\end{eqnarray*}
Then each set $\cu_h \setminus \cs_h^j$ is connected.
A component of $\cs_h^j$ is just the saturation of a
component of $\cc_h^j$ which is not enclosed by any component
of $\bigcup_{g=h+1}^l \cc_g^j$.
Each component $\ck$ of $\cs_h^j$ has piecewise smooth boundary, and
according to Lemma~\ref{l:nu}\:(ii) it contains only one $h$-cube
$\psi_h(C)$, and
\begin{equation*}  
\psi_h^{-1} \left( \ck \right) \,\subset\, \cn_{\nu_h d_h} (C),
\quad \, h=1, \dots, l.
\end{equation*}
While a component of $\cs_0^j$ is a cube of $\cc_0^j$ and a component
of $\cs_1^j$ is the union of a cube of $\cc_1^j$ and the overlapping
cubes of $\cc_0^j$, 
a component of $\cs_2^j$ might contain cubes of $\cc_0^j$ which are
disjoint from $\cc_1^j \cup \cc_2^j$, cf.~Figure~\ref{figure28}.

%
%

\ni
If the ratios $d_h/d_{h+1}$, $h=0,\dots,l-1$, are small enough, then
Lemma~\ref{l:nu}\:(ii) implies that a component of $\cc_h^j$ cannot be
enclosed by a component of $\cc_g^j$ for some $g<h$, and so the sets
$\cs_h^j$, $h=0,\dots,l$, are disjoint.
We finally abbreviate
\[
\cs^j \,:\,=\, \bigcup_{h=0}^l \cs_h^j
\]
and read off from~\eqref{id:MUU} and the definition of the sets
$\cs_h^j$ that
\begin{eqnarray}  \label{id:MSj} 
M \,=\, \bigcup_{j=1}^k \cs^j .
\end{eqnarray}

\b
\ni
{\bf Step~4. Moving the cubes of the same colour into $\cb_0$} 

\begin{figure}[ht] 
 \begin{center}
  \leavevmode\epsfbox{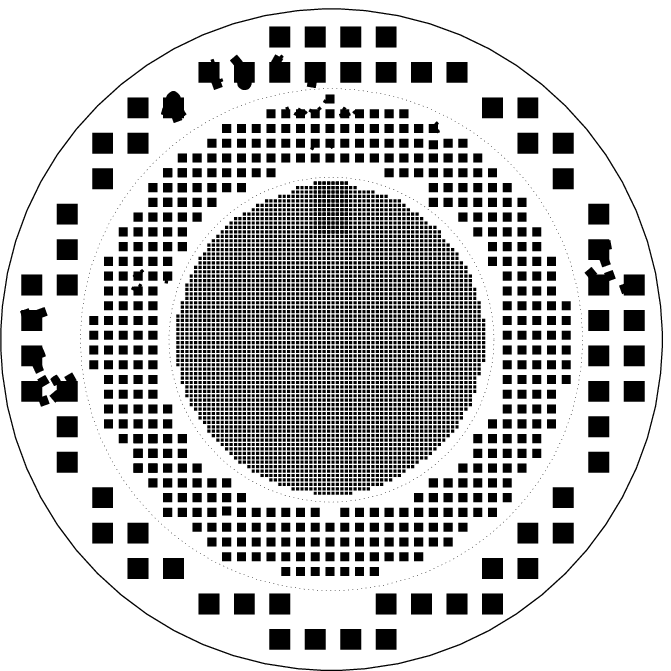}
 \end{center}
 \caption{The image $\left( \psi_0^{-1} \circ \Phi^j\right) (\cs^j) \subset 
\psi_0^{-1}(\cb_0)$ for $l=2$.} 
 \label{figure.idea}
\end{figure}

\s
\ni
In order to move the sets $\cs^j$ into $\cb_0$
we shall possibly have to choose the $d_i$'s yet smaller.
We shall then be able to construct for each $j$ a Hamiltonian isotopy $\Phi^j$ 
of $M$
which first moves $\cs_0^j$ to a ``dense cluster'' around the center
of $\cb_0$ and then successively moves $\cs_h^j$ to a ``shell'' around
the already constructed cluster 
$\bigcup_{g=0}^{h-1} \Phi^j \bigl( \cs_g^j \bigr)$, 
$h = 1, \dots, l$, 
cf.~Figure~\ref{figure.idea}.

%
%

\ni
The main tool for the construction of the maps $\Phi^j$ is the
following elementary lemma.

\begin{lemma}  \label{l:trans}
Let $K$ be a compact subset of $\RR^{2n}$ and let $q$ be a point in
$\RR^{2n}$. 
Denote by $\ck$ the convex hull of the union $K \cup (q + K)$.
For any open neighbourhood $U$ of $\ck$ there exists a
symplectomorphism $\tau$ of $\RR^{2n}$ which is supported in $U$ and
which translates $K$ to $q+K$.
\end{lemma}

\proof
We follow \cite[p.~73]{HZ}.
We choose a smooth function $f \colon \RR^{2n} \ra \RR$ such that $f
|_{\ck} =1$ and $f |_{\RR^{2n} \setminus U} =0$.
Define the Hamiltonian function $H \colon \RR^{2n} \ra \RR$ by
\[
H(z) \,=\, f(z) \langle z, -Jq \rangle 
\]
where $\langle \cdot, \cdot \rangle$ denotes the Euclidean scalar
product on $\RR^{2n}$ and where $J$ denotes the standard complex
structure on $\RR^{2n}$ defined by
\[
\omega_0 (z,w) \,=\, \langle z, -Jw \rangle, \quad \, z,w \in
\RR^{2n} . 
\]
Recall that the Hamiltonian vector field $X_H$ of $H$ is given by
$X_H (z) = J \nabla H (z)$.
We conclude that the time-1-map $\tau$ of the flow generated by $X_H$
is a symplectomorphism of $\RR^{2n}$ which is supported in $U$.
Moreover, for $z \in \ck$ we have
\[
X_H (z) \,=\, J \nabla H (z) \,=\, J (-Jq) \,=\, q, 
\]
and so $\tau (z) = z+q$ for all $z \in K$.
\proofend

We denote by $B_r$ the open ball in $\RR^{2n}$ of radius $r$ and centered 
at the origin.
We recursively define the open balls $B_{r_0}, \ldots, B_{r_l}$ and the open ``annuli'' 
$A_{r_{h-1}}^{r_h} = B_{r_h} \setminus \overline{B_{r_{h-1}}}$ by
\begin{eqnarray}
\left| B_{r_0} \right|           &=& \tfrac{1}{k} \left( \left| V_0 \right| +
    \tfrac{k-1}{l+1} \eps \right) , \label{def:r0}\\
\left| A_{r_{h-1}}^{r_h} \right| &=& \tfrac{1}{k} \left( \left|V_h \right| +
    \tfrac{k-1}{l+1} \eps \right) , \quad \; h=1, \dots, l. \label{def:ri}
\end{eqnarray}
The definitions \eqref{def:r0} and \eqref{def:ri}, the
identities $\left| V_h \right| = \mu \left( \cu_h \right)$ and the
  estimate \eqref{est:mme}, and the estimate \eqref{est:B0p} and the identity
  $\left| B^{2n}(a_0') \right| = \mu \left( \cb_0' \right)$ imply that
\begin{eqnarray}  \label{est:vol} 
 \left| B_{r_0} \right| + \sum_{h=1}^l \left| A_{r_{h-1}}^{r_h} \right|
  &=& \frac{1}{k} \sum_{h=0}^l  \left( \left| V_h \right| + 
                                      \frac{k-1}{l+1} \eps \right)  \\
  &<& \frac{\mu (M)}{k} + \frac{\eps}{k} + \frac{k-1}{k} \eps \notag\\
  &=& \frac{\mu (M)}{k} + \eps \notag\\
  &<& \left| B^{2n}(a_0') \right| \notag
\end{eqnarray}
and so 
\begin{equation}  \label{inc:BAB}
B_{r_0} \cup \bigcup_{h=1}^l A_{r_{h-1}}^{r_h} \,\subset\, B^{2n}(a_0').
\end{equation}

Consider again the symplectic atlas
$\ea = \left\{ \psi_h \colon V_h \ra \cu_h, \: h = 0, \dots, l
\right\}$ of $(M, \go)$ constructed in Step~1. 
Recall that $\psi_0 \colon V_0 \ra \cu_0$ is
the Darboux chart $\gf_0 \colon B^{2n}(a_0) \ra \cb_0$ and that the
sets $\cu_h$ and $V_h$ are connected and have piecewise smooth
boundaries.
Also recall that there exist larger charts $\widetilde{\psi}_h \colon
\widetilde{V}_h \ra \widetilde{\cu}_h$.
We can assume that the sets $\widetilde{\cu}_h$ and $\widetilde{V}_h$
are also connected and have piecewise smooth boundaries. 
We fix $j \in \left\{ 1, \dots, k \right\}$.
The construction of the map $\Phi_0^j$ will somewhat differ from the
one of the maps $\Phi_h^j$ for $h \ge 1$ since $\Phi_0^j \bigl( \cs_0^j
\bigr)$ will not be disjoint from $\cs_0^j$. We start with
constructing $\Phi_0^j$.

\begin{proposition}  \label{p:Phi0}
If the numbers $d_0, \dots, d_l >0$ 
are small enough,
then there exists a symplectomorphism $\Phi_0^j$ of $M$ whose support
is disjoint from $\bigcup_{h=1}^l \cs_h^j$ and such that
$\Phi_0^j \bigl( \cs_0^j \bigr) \subset \psi_0 \left( B_{r_0} \right)$. 
\end{proposition} 

\proof
We recall that $\cs_0^j$ is the union of ``free'' cubes of $\cc_0^j$, i.e.,
each component of $\cs_0^j$ is a cube of $\cc_0^j$ 
which is not enclosed by any component of $\bigcup_{h=1}^l \cc_h^j$.
We abbreviate $\es_0 = \psi_0^{-1} \bigl( \cs_0^j \bigr)$. 
So, $\es_0$ is a disjoint union of cubes. 
Since $\es_0$ is contained in $\ec_0^j (d_0)$,
we deduce from the estimate \eqref{e:Cijdi} for $i=0$ and from
definition~\eqref{def:r0} that
\begin{equation}  \label{est:S0jB}
\left| \es_0 \right| \,<\, \left| B_{r_0} \right| .
\end{equation}
We denote by $\eq$ the standard decomposition of $\RR^{2n}$ into closed cubes,
\[
\eq \,:\,=\, \left\{v+ [0,1]^{2n}\;\Big\vert\;  v \in \ZZ^{2n}\right\}.
\] 
Furthermore, for each $\nu>0$ we set 
\[
\nu\:\!\eq \,:\,=\, \left\{\nu \:\!v+ [0,\nu]^{2n}\;\Big\vert \;v \in \ZZ^{2n}\right\},
\] 
and for each subset $A$ of $\RR^{2n}$ we denote by
$\eq (\nu, A)$ the union of cubes of $\nu \:\!\eq$ which are contained in $A$.
By ``a cube of $\eq (\nu, A)$" we shall mean a cube of $\nu \:\!\eq$ contained in $A$.
Let $s_0$ be the number of components (i.e., cubes) of $\es_0$.
The estimate~\eqref{est:S0jB} implies that after choosing $d_0 >0$
smaller if necessary we find $\eps_0>0$ such that  
$\eq \left(d_0+\eps_0, B_{r_0} \right)$
contains at least $s_0$ cubes.

Recall that $k \ge 2n+1$ and recall from the estimate \eqref{est:vol}
that $r_0 < \sqrt{a_0'/\pi}$. We define $\widetilde{r}_0 > r_0$ by
\begin{equation}  \label{def:r0+}
\widetilde{r}_0 \,=\, \min \left\{ \frac{2k}{4n+1} r_0,\, \frac{1}{2} \left(
                                   r_0 + \sqrt{a_0'/\pi} \right) \right\} 
\end{equation}
and we denote by $\es_0^{\inter}$ the
union of those cubes of $\es_0$ which are contained in $B_{\widetilde{r}_0}$.
Since $B_{\widetilde{r}_0} \subset B^{2n}(a_0')$ and since $\cb_0' = \psi_0
\left( B^{2n}(a_0') \right)$ is disjoint from $\cu_h$ and $\cs_h^j
\subset \cu_h$, $h \ge 1$, the set $B_{\widetilde{r}_0}$ is disjoint from
$\psi_0^{-1} \bigl( \cs_h^j \bigr)$, $h \ge 1$. 
In particular, $\es_0^{\inter}$ is the union of cubes of 
$\ec_0^j(d_0)$ contained in $B_{\widetilde{r}_0}$, 
cf.~Figure~\ref{figure32}.
We abbreviate the union of exterior cubes of $\es_0$ by
\[
\es_0^{\ext} \,:\,=\, \es_0 \setminus \es_0^{\inter} .
\]

\begin{lemma}  \label{l:theta}
For $d_0$ and $\eps_0$ small enough there exists a symplectomorphism
$\theta$ of $\widetilde{V}_0$ such that 
\begin{itemize}
\item[(i)]
the support of $\theta$ is contained in $B_{\widetilde{r}_0}$ and disjoint from
$\es_0^{\ext}$; 
\vspace{.2em}
\item[(ii)]
$\theta$ maps each cube of $\es_0^{\inter}$ into a cube of 
$\eq \left( d_0+\eps_0, B_{r_0} \right)$;
\vspace{.4em}
\item[(iii)]
the union of cubes of $\eq \left( d_0+\eps_0, B_{r_0}
\right)$ containing a cube of $\theta \left( \es_0^{\inter} \right)$
is contractible.
\end{itemize}
\end{lemma}

\begin{figure}[ht] 
 \begin{center}
  \psfrag{1}{$\theta_1$}
  \psfrag{2}{$\theta_2$}
  \psfrag{B}{$B_{r_0}$}
  \psfrag{Bt}{$B_{\widetilde{r}_0}$}
  \leavevmode\epsfbox{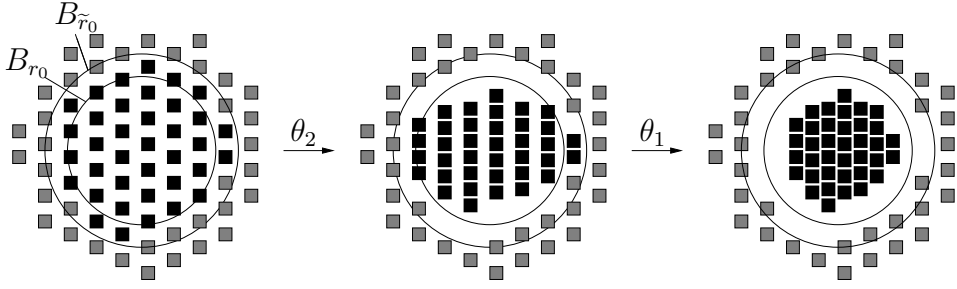}
 \end{center}
 \caption{The map $\theta = \theta_1 \circ \theta_2$ for $j=1$.} 
 \label{figure30}
\end{figure}

\begin{figure}[ht] 
 \begin{center}
  \psfrag{r0t}{$\widetilde{r}_0$}
  \psfrag{a'}{$\sqrt{\frac{a_0'}{\pi}}$}
  \psfrag{a}{$\sqrt{\frac{a_0}{\pi}}$}
  \leavevmode\epsfbox{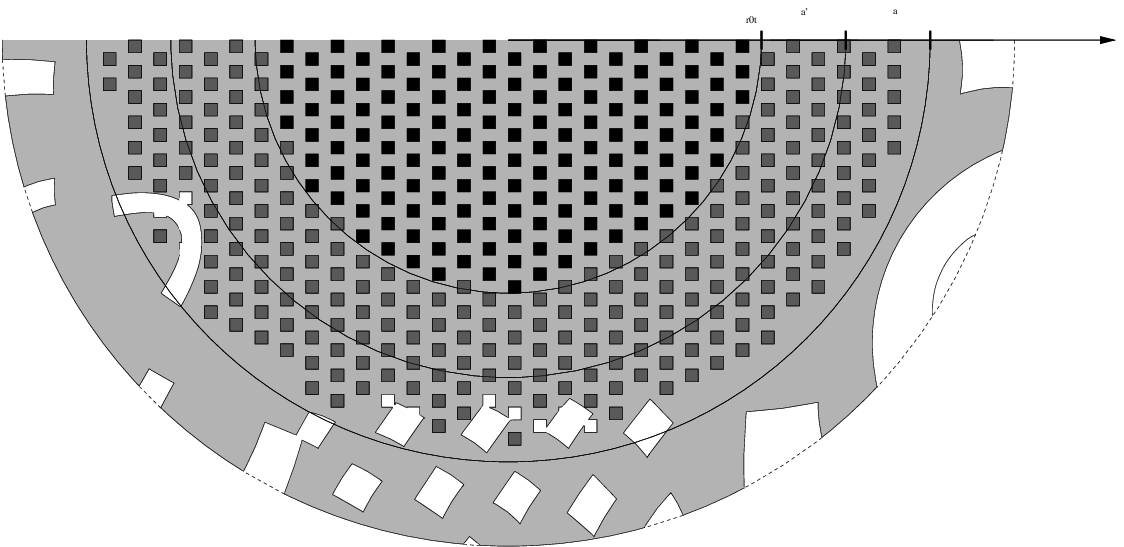}
 \end{center}
 \caption{Half of the subset $\es_0 = \es_0^{\inter} \cup \es_0^{\ext}$ 
                                             of $\widehat{V}_0$.}  
 \label{figure32}
\end{figure}

\proof
Using Lemmata~\ref{l:cover} and \ref{l:trans} we successively construct
symplectomorphisms $\theta_{2n}, \theta_{2n-1}, \dots, \theta_1$
such that $\theta_{2n}$ ``compresses'' $\es_0^{\inter}$ along the
$x_{2n}$-axis and $\theta_i$ ``compresses'' $\theta_{i+1} \circ \dots
\circ \theta_{2n} \left( \es_0^{\inter} \right)$ along the
$x_i$-axis, $i = 2n-1, \dots, 1$, and such that the composite map
\[
\theta \,=\, \theta_1 \circ \dots \circ \theta_{2n}
\]
meets assertion (i) as well as assertions (ii) and (iii) with $\eq
\left( d_0 + \eps_0, B_{r_0} \right)$ replaced by $\eq \left(
  d_0+\eps_0, B_{\widetilde{r}_0} \right)$, 
cf.~Figure~\ref{figure30}.

%
%

\ni
In order to see that assertions (ii) and (iii) can be fulfilled as
stated, we infer from the definition of the set $d_0 \ec^j (2n,k)
\supset \es_0^{\inter}$ given in Step~2 that
\[
\frac{\diameter \es_0^{\inter}}{\diameter \theta \left(
    \es_0^{\inter} \right) } 
\,\ra\, \frac{k}{2n} \quad \text{ as }\: d_0 \ra 0 \text{ and } \eps_0
\ra 0.
\]
In view of the choice \eqref{def:r0+} of $\widetilde{r}_0$ we can
therefore choose $d_0$ and $\eps_0$ so small that 
$\theta \left( \es_0^{\inter} \right) \subset 
\eq \left( d_0+\eps_0, B_{r_0} \right)$, as desired. 
\proofend

\begin{lemma}  \label{l:thetaj}
If the numbers $d_0, \dots, d_l >0$ 
are small enough, then
there exists a symplectomorphism
$\Theta_0$ of $\widetilde{V}_0$ such that 
\begin{itemize}
\item[(i)]
the support of $\Theta_0$ is compact and disjoint from 
\[
\psi_0^{-1} \left( \bigcup_{h=1}^l \cs_h^j \right) \, \cup \; \theta
\left( \es_0^{\inter} \right) ;
\] 
\vspace{.2em}
\item[(ii)]
$\Theta_0$ maps each cube of $\es_0^{\ext}$
into a cube of $\eq \left( d_0+\eps_0, B_{r_0} \right)$.
\end{itemize}
\end{lemma}

\proof
The set $\cu_0 \setminus \bigcup_{h=1}^l \cs_h^j$ might not be
connected for any choice of $d_0, \dots, d_l$, in which case not every
cube of $\cs_0^j$ can be moved into $\psi_0 \left( B_{r_0} \right)$
inside $\cu_0 \setminus \bigcup_{h=1}^l \cs_h^j$,
cf.~Figure~\ref{figure32}.
This is the reason
why we work in the extended chart $\widetilde{\psi}_0 \colon
\widetilde{V}_0 \ra \widetilde{\cu}_0$. 
We choose the numbers $d_0, \dots, d_l$ so small that each component
of $\bigcup_{h=1}^l \cs_h^j$ which intersects $\cu_0$ is contained in
$\widetilde{\cu}_0$. The component $\widehat{\cu}_0$ of
$\widetilde{\cu}_0 \setminus \bigcup_{h=1}^l \cs_h^j$ containing
$\cb_0'$ then contains $\cs_0^j$, and the set 
$\widehat{V}_0 := \widetilde{\psi}_0^{-1} \bigl( \widehat{\cu}_0 \bigr)$ is an 
open
connected set with piecewise smooth boundary which contains $\es_0$,
cf.~Figure~\ref{figure32}.
(In this figure, one can still find one component of $\bigcup_{h=1}^l \cs_h^j$ 
which intersects $\cu_0$ but is not contained in $\widetilde{\cu}_0$!)

%
%

\begin{figure}[ht] 
 \begin{center}
  \psfrag{1}{$1$}
  \psfrag{2}{$2$}
  \leavevmode\epsfbox{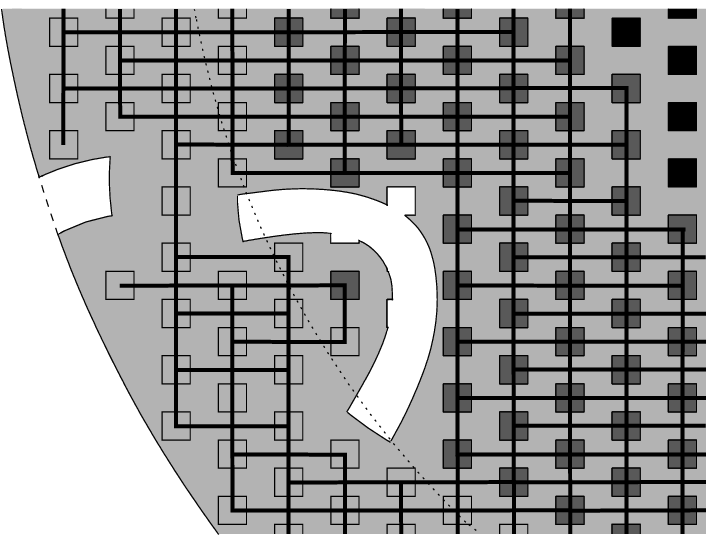}
 \end{center}
 \caption{Part of the graph $\cg_0$ associated with 
          $\widehat{\es}_0^{\ext}$.} 
 \label{figure31}
\end{figure}

\ni
In order to move the cubes of $\es_0^{\ext}$ into $B_{r_0}$ we shall
associate a tree with $\es_0^{\ext}$.
Recall that $\es_0^{\ext}$ is a subset of $d_0 \ec^j (2n,k)$. We
enlarge $\es_0^{\ext}$ to the set $\widehat{\es}_0^{\ext}$
defined as the union of cubes of $d_0 \ec^j(2n,k) \setminus
\es_0^{\inter}$ which are contained in $\widehat{V}_0$. 
Abbreviate
\[
\gl_m \,:\,=\,                     
    \left\{
         \begin{array}{ll}   
                k      & \text{if }\, m = 1,  \\
                2n-m+2 & \text{if }\, m \in \left\{ 2, \dots, 2n \right\} .
         \end{array}   \right. 
\]  
We say that two cubes $C$ and $C'$ of $\widehat{\es}_0^{\ext}$
are {\it $m$-neighbours}\, if 
\[
C' \,=\, C \pm d_0 \gl_m e_m
\]
for some $m \in \left\{ 1, \dots, 2n \right\}$ and if the convex
hull of $C \cup C'$
is contained in $\widehat{V}_0$. According to Lemma~\ref{l:cover}\:(i) the   
interior of the convex hull of two $m$-neighbours does not intersect any
third cube of $\widehat{\es}_0^{\ext}$, cf.~Figure~\ref{figure20}.
We define $\cg_0'$ to be the graph whose edges are the straight segments 
joining the centers of neighbours in $\widehat{\es}_0^{\ext}$,
and we define $\cg_0$ to be the graph obtained from $\cg_0'$ by
declaring the intersections of edges to be vertices, 
cf.~Figure~\ref{figure31}. 

%
%

\ni
Since $\widehat{V}_0$ is an open connected relatively compact set
with piecewise smooth boundary, we can choose $d_0$ so small that 
the graph $\cg_0$ is connected.
Choosing $d_0$ yet smaller if necessary, we can also assume that
\begin{equation}  \label{est:drr}
\sqrt{2n} \, d_0 \,<\, \frac{\widetilde{r}_0-r_0}{2} 
\end{equation}
and that the convex hull of the union $C \cup C'$ of any two neighbours in 
$\widehat{\es}_0^{\ext}$ is 
contained in $\widehat{V}_0 \setminus \overline{B_{r_0}}$.
We then in particular have that
$\es_0^{\ext}$ is disjoint from $\overline{B_{r_0}}$.
Let $C_1$ be a cube of $\es_0^{\ext}$ whose distance to $B_{r_0}$ is
minimal. 
We choose a maximal tree $\ct_0$ in $\cg_0$ which is rooted at the
center of $C_1$. Denote a vertex of $\ct_0$ represented by the center
of a cube $C$ of $\es_0^{\ext}$ by $v(C)$
and write $\prec$ for the partial ordering on $\es_0^{\ext}$ induced
by $\ct_0$.
We number the $s_0^{\ext}$ many cubes of $\es_0^{\ext}$ in such a
way that
\begin{equation}  \label{inc:v}
v(C_c) \prec v(C_{c'}) \;\Longrightarrow\; c<c' .
\end{equation}
We finally recall that $\eq \left( d_0+\eps_0, B_{r_0} \right)$
contains at least $s_0$ cubes.
Denote by $\eq \left( \theta \left( \es_0^{\inter} \right) \right)$ 
the union of those cubes in 
$\eq(d_0+\eps_0, B_{r_0})$ which contain a cube of 
$\theta \left( \es_0^{\inter} \right)$.
According to Lemma~\ref{l:theta}\:(iii), the set
$\eq \left( \theta \left( \es_0^{\inter} \right) \right)$
is contractible.
We can therefore successively choose cubes 
$Q_1, \dots, Q_{s_0^{\ext}}$ from $\eq(d_0+\eps_0, B_{r_0})$
different from the cubes of 
$\eq \left( \theta \left( \es_0^{\inter} \right) \right)$
in such a way that each of the sets
\begin{equation}  \label{s:scon}
\eq \left( \theta \left( \es_0^{\inter} \right) \right)
\,\cup\, \bigcup_{b=1}^c Q_b,       
\end{equation}
$c=1,\dots,s_0^{\ext}$,
is contractible.

We are now in a position to move the cubes of $\es_0^{\ext}$ into $B_{r_0}$.
We shall successively move $C_c$ into $Q_c$, 
$c=1, \dots, s_0^{\ext}$. Define $\widehat{r}_0 \in \;]r_0,
\widetilde{r}_0[$ by
$\widehat{r}_0 := (r_0+\widetilde{r}_0)/2$.
In view of assumption~\eqref{est:drr} we can then estimate the
diameter of a cube of $\es_0^{\ext}$ by
\begin{equation}  \label{est:d0rr}
\sqrt{2n}\, d_0 \,<\, \widehat{r}_0 - r_0 .
\end{equation}
We first use Lemma~\ref{l:trans} to construct a symplectomorphism
$\gt_1$ of $\widetilde{V}_0$ whose support is contained in
$\widehat{V}_0$ and is disjoint from  
\[
\bigcup_{b=2}^{s_{0}^{\ext}}C_b \,\cup\, \theta \left(
  \es_0^{\inter} \right)
\] 
and which maps $C_1$ into $Q_1$.
Indeed, since $C_1$ is a cube of $\es_0^{\ext}$ closest to $B_{r_0}$
and in view of the estimate~\eqref{est:d0rr}, 
we can first move $C_1$ into the annulus
$B_{\widehat{r}_0} \setminus B_{r_0}$ without touching
$\bigcup_{b\ge2} C_b$, and since $\eq \left( \theta \left( \es_0^{\inter} \right) \right)$
is contractible,
we can then move the image cube along a piecewise linear path inside 
$B_{\widehat{r}_0} \setminus B_{r_0}$ to a position from which it can
be moved into $B_{r_0}$ to its preassigned cube $Q_1$ without
touching $\theta \left( \es_0^{\inter} \right)$.

\begin{figure}[ht]                     
 \begin{center}
  \psfrag{1}{$C_-$}
  \psfrag{2}{$C_+$}
  \psfrag{3}{$C_a$}
  \psfrag{4}{$C_{a'}$}
  \leavevmode\epsfbox{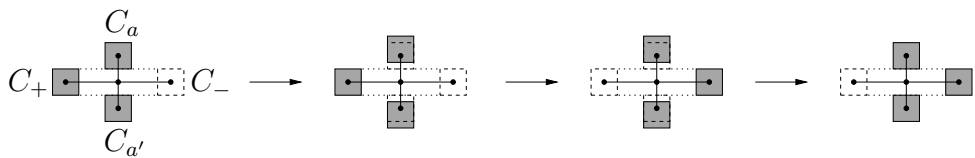}
 \end{center}
 \caption{How to move $C_+$ to $C_-$ along a path blocked by $C_a$ and
   $C_{a'}$.} 
 \label{figure33}
\end{figure}

Assume now by induction that we have already constructed
symplectomorphisms $\gt_b$ which moved the cubes $C_b$ into
the cubes $Q_b$ for $b=1,\dots,c-1$. We are going to construct a
symplectomorphism $\gt_c$ of $\widetilde{V}_0$ whose support
is contained in $\widehat{V}_0$ and is disjoint from 
\begin{equation}  \label{s:CQt}
\bigcup_{b=c+1}^{s_0^{\ext}} C_b \,\cup\, \bigcup_{b=1}^{c-1} Q_b
\,\cup\, \theta \left( \es_0^{\inter} \right)
\end{equation}
and which maps $C_c$ into $Q_c$.
Let $\gg$ be the piecewise linear path from $v(C_c)$ to $v(C_1)$
determined by the tree $\ct_0$.
Because of~\eqref{inc:v}, all the cubes of $\es_0^{\ext}$ on $\gg$
except $C_c$ have already been moved into $B_{r_0}$. 
Using Lemmata~\ref{l:cover}\:(i) and \ref{l:trans} we can therefore move $C_c$ 
along
$\gg$ to (the ``former locus'' of) $C_1$ without touching 
$\bigcup_{b \ge c+1} C_b$.
More precisely, consider two consecutive cubes $C_+$ and $C_-$ along $\gg$ which are centred at the vertices $+$ and $-$ of $\cg_0$, respectively, and let $\gs$ be the part of $\gg$ joining $+$ and $-$.
As the notation suggests, $- \prec +$ along $\gg$.
The path $\gs$ may consist of one or two or more than two edges, which may be parallel to different coordinate axes, cf.~Figure~\ref{figure21}.
We only describe a typical case, in which $\gs$ consists of two edges parallel to the same coordinate axis.
Let $R$ be the convex hull of $C_+ \cup C_-$.
In view of Lemma~\ref{l:cover}\:(i), the closed rectangle $R$ either is disjoint 
from
$\bigcup_{b \ge c+1} C_b$ or it touches some cubes $C_a$ with $a \ge c+1$
along a face.
In the first case, we can directly apply Lemma~\ref{l:trans} to move
$C_+$ to $C_-$ without touching $\bigcup_{b \ge c+1} C_b$. In the
second case, we first move the touching cubes $C_a$ a bit away from
$R$, then move $C_+$ to $C_-$, and then move the displaced cubes
back to their former locus, cf.~Figure~\ref{figure33}. 

%
%

\ni
We can do this in such a way that the support of the resulting map 
$\tau_{\gs}$ which translates $C_+$ to $C_-$ is disjoint from 
$\bigcup_{b \ge c+1} C_b$.
Since $R$ is contained in $\widehat{V}_0 \setminus \overline{B_{r_0}}$
we can also arrange the support of $\tau_{\gs}$ to be contained in 
$\widehat{V}_0 \setminus \overline{B_{r_0}}$. 
Composing the maps $\tau_{\gs}$ corresponding to the parts $\gs$ of
$\gg$ we obtain a symplectomorphism $\tau_c$ whose support is
contained in $\widehat{V}_0$ and is disjoint from the set~\eqref{s:CQt}
and which maps $C_c$ to $C_1$.
Since the set~\eqref{s:scon} is contractible, we can now
proceed as in the construction of $\gt_1$ and construct a
symplectomorphism $\gt_c$ which moves the image
of $C_c$ at $C_1$ into $Q_c$ without touching the set~\eqref{s:CQt}.
The composition $\gt_c \circ \tau_c$ is as desired.

After all, the composite map
\[
\Theta_0 \,=\, \left( \gt_{s_0^{\ext}} \circ \tau_{s_0^{\ext}} \right)
\circ \dots \circ
\left( \gt_2 \circ \tau_2 \right) \circ \gt_1
\]
is a symplectomorphism of $\widetilde{V}_0$ which meets assertions~(i)
and (ii).
\proofend

Let $\theta$ and $\Theta_0$ be the symplectomorphisms guaranteed
by Lemmata~\ref{l:theta} and \ref{l:thetaj}.
The symplectomorphism
\[
\widetilde{\psi}_0 \circ \Theta_0 \circ \theta \circ \widetilde{\psi}_0^{-1}
\]
of $\widetilde{\cu}_0$ smoothly extends by the identity to a
symplectomorphism $\Phi_0^j$ of $M$ whose support
is disjoint from $\bigcup_{h=1}^l \cs_h^j$ and such that
$\Phi_0^j \bigl( \cs_0^j \bigr) \subset \psi_0 \left( B_{r_0} \right)$.
The proof of Proposition~\ref{p:Phi0} is complete.
\proofend

\begin{proposition}  \label{p:Phih}
If the numbers $d_0, \dots, d_l >0$ as well as the ratios 
$d_i / d_{i+1}$, $i = 0, \dots, l-1$, are small enough,
then for each $h = 1, \dots,l$ there exists a symplectomorphism
$\Phi_h^j$ of $M$ whose support is disjoint from
\[
\bigcup_{g=0}^{h-1} \Phi_g^j \bigl( \cs_g^j \bigr) \cup
\bigcup_{g=h+1}^l \cs_g^j
\]
and such that
$\Phi_h^j \bigl( \cs_h^j \bigr) \subset \psi_0 \bigl( A_{r_{h-1}}^{r_h} \bigr)$. 
\end{proposition}

\begin{figure}[ht] 
 \begin{center}
 \leavevmode\epsfbox{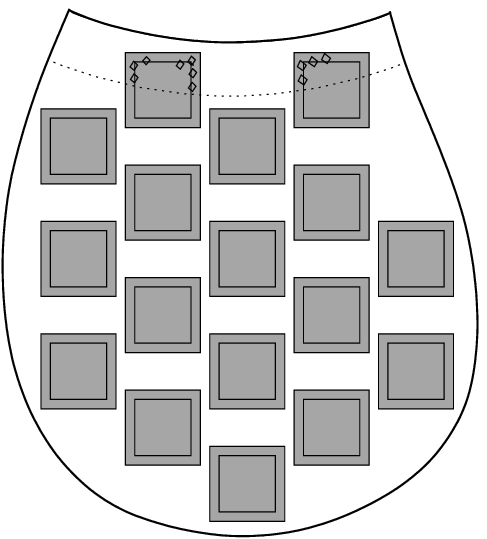}
\end{center}
\caption{The set $\es_1 \subset V_1$.} 
 \label{figuresigma1}
\end{figure}

\proof
We first explain the construction of $\Phi_1^j$.
Recall from the end of Step~3 that $\cs_1^j \subset \cu_1$ is the
union of those components of $\cc_1^j$ which are not enclosed by any
component of $\bigcup_{h=2}^l \cc_h^j$.
Each component $\ck$ consists of a $1$-cube $\psi_1(C)$ and some
overlapping cubes of $\cc_0^j$, and
\[
\psi_1^{-1} (\ck) \,\subset\, \cn_{\nu_1 d_1} (C) \,\subset\, V_1 = 
\psi_1^{-1}(\cu_1)
\]
according to~\eqref{inc:NV} and Lemma~\ref{l:nu}\:(ii).
For any cube $C$ of $d_1 \ec^j (2n,k)$ we denote by $C^{\nu_1}$ the
closed cube of width $(1+2\nu_1) d_1$ concentric to $C$. 
If $C$ belongs to $\ec_1^j (d_1)$, then $C^{\nu_1}$ is 
the
smallest closed cube containing the neighbourhood $\cn_{\nu_1 d_1}
(C)$ of $C$. We abbreviate
\[
\es_1 \,:\,=\, \bigcup C^{\nu_1}
\]
where the union is taken over those cubes $C$ of $\ec_1^j(d_1)$ that
lie in $\psi_1^{-1} \bigl( \cs_1^j \bigr)$, 
see Figure~\ref{figuresigma1}.

%
%

\ni
In view of the choice~\eqref{est:nud} the cubes $C^{\nu_1}$ are disjoint. 
Since the compact subset $\psi_1^{-1} \bigl( \cs_1^j \bigr)$ of $V_1$ is
disjoint from the compact subset 
$\psi_1^{-1} \bigl( \bigcup_{h=2}^l \cs_h^j \bigr)$ of $\overline{V_1}$,
we can choose $\nu_1 >0$ (and for this $d_0>0$) so small that
$\es_1$ is disjoint from $\psi_1^{-1} \bigl( \bigcup_{h=2}^l \cs_h^j \bigr)$.
Since for each cube $C^{\nu_1}$ of $\es_1$ the cube $C$ belongs to 
$\ec_1^j(d_1)$, we read off from estimate~\eqref{est:2nu} for
$i=1$ and from definition~\eqref{def:ri} for $h=1$ that 
\begin{equation}  \label{est:S1A}
\left| \es_1 \right| \,<\, \left| A_{r_0}^{r_1} \right| .
\end{equation} 
Let $s_1$ be the number of cubes of $\es_1$. The estimate~\eqref{est:S1A} implies that after choosing $d_1>0$ and $\nu_1 >0$
smaller if necessary we find $\eps_1>0$ such that 
$\eq \left( (1+2\nu_1)d_1 + \eps_1, A_{r_0}^{r_1} \right)$ 
contains at least $s_1$ cubes.

\begin{figure}[ht] 
 \begin{center}
  \psfrag{pB}{$\psi_0 (B_{r_l})$}
  \psfrag{U0'}{$\cu_0'$}
  \psfrag{U0}{$\cu_0$}
  \psfrag{W0}{$\cw_0$}
  \psfrag{p1}{$p_1$}
  \psfrag{U1'}{$\cu_1'$}
  \psfrag{U1}{$\cu_1$}
  \psfrag{U1t}{$\widetilde \cu_1$}
 \leavevmode\epsfbox{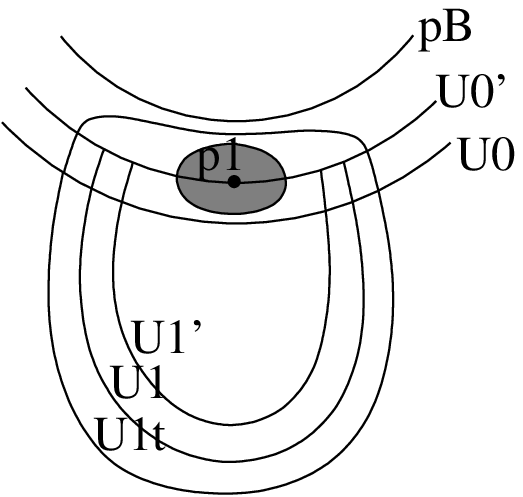}
\end{center}
\caption{The neighbourhood $\cw_1$ of $p_1$.} 
 \label{figureW1}
\end{figure}

Recall from Step~1 that $\psi_1 \colon V_1 \to \cu_1$ is the restriction of a 
larger chart $\widetilde \psi_1 \colon \widetilde V_1 \to \widetilde \cu_1$.
In view of~\eqref{est:vol} we have 
\begin{equation}  \label{eq:rl}
r_l < \sqrt{a_0' / \pi} ,
\end{equation}
and so we can assume that $\widetilde \cu_1$ is disjoint from 
$\psi_0 (B_{r_l})$.
Also recall that there exists a point $p_1 \in \pp \cu_0' \cap \pp \cu_1'$ and a 
neighbourhood $\cw_1 \subset \widetilde \cu_0 \cap \widetilde \cu_1$ of $p_1$ 
such that 
$\widetilde\psi_0^{-1} \circ \widetilde \psi_1$ restricts to the identity on 
$\cw_1$, see Figure~\ref{figureW1}.

%
%

\ni
Since $\cs_h^j$ is a compact subset of $\cu_h$, and since $\cu_h$ is 
disjoint from $\cu_0'$ and $\cu_1'$ for $h \ge 2$ 
according to~\eqref{eq:uiuj}, 
the point $p_1$ is disjoint from $\bigcup_{h=2}^l \cs_h^j$.
We choose the numbers $d_0, \dots, d_l$ so small that each component
of $\bigcup_{h=2}^l \cs_h^j$ which intersects $\cu_1$ is contained in
$\widetilde{\cu}_1$. The component $\widehat{\cu}_1$ of
$\widetilde{\cu}_1 \setminus \bigcup_{h=2}^l \cs_h^j$ containing $p_1$ 
then contains $\cs_1^j$, and the set 
$\widehat{V}_1 := \widetilde{\psi}_1^{-1} \bigl( \widehat{\cu}_1 \bigr)$ is an 
open
connected set with piecewise smooth boundary which contains $\es_1$.
After choosing $\cw_1$ smaller if necessary, we can assume that
$\cw_1 \subset \widehat \cu_1$ and 
$W_1 := \widetilde \psi_1^{-1}(\cw_1) \subset \widehat V_1$.

\begin{figure}[ht] 
 \begin{center}
  \psfrag{WV}{$W_1 \cap V_1$}
  \psfrag{Wp}{$W_1 \cap \widetilde \psi_1^{-1}(\cu_0')$}
 \leavevmode\epsfbox{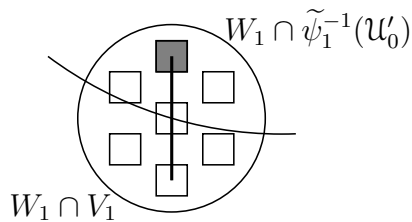}
\end{center}
\caption{The pilot cube $C_{\frak p}^{\nu_1} \subset W_1$ and 
         the linear tree $\ct_1'$.} 
 \label{figureCp}
\end{figure}

We enlarge $\es_1$ to the set $\widehat{\es}_1 := \bigcup C^{\nu_1}$
where the union is taken over all cubes $C$ of $d_1 \ec^j(2n,k)$ 
which are contained in $\widehat{V}_1$.
In the same way as in the proof of Lemma~\ref{l:thetaj} we associate a graph 
$\cg_1$ to $\widehat{\es}_1$, which is connected for $d_0, d_1$ small enough.
Choosing $d_0, d_1$ yet smaller if necessary, 
we find a linear tree $\ct_1' \subset \cg_1$
which is contained in $W_1$,
is rooted in the center of a ``pilot cube'' 
$C_{\frak p}^{\nu_1} \subset \widetilde \psi_1^{-1}(\cu_0')$, 
and meets at least one cube of $\es_1$,
see Figure~\ref{figureCp}.

%
%

\ni
Choose a maximal tree $\ct_1 \subset \cg_1$ which is also rooted in 
$C_{\frak p}^{\nu_1}$ and contains $\ct_1'$.
Denoting a vertex of $\ct_1$ represented by the center
of a cube $C^{\nu_1}$ of $\es_1$ by $v(C^{\nu_1})$
and writing $\prec$ for the partial ordering on $\es_1$ induced by $\ct_1$,
we number the $s_1$ many cubes of $\es_1$ in such a way that
\begin{equation*}  
v(C_c^{\nu_1}) \prec v(C_{c'}^{\nu_1}) \;\Longrightarrow\; c<c' .
\end{equation*}

We are now in a position to move the $s_1$ cubes of $\es_1$ into 
$A_{r_0}^{r_1}$. Assume by induction that we have already constructed 
symplectomorphisms $\theta_b$, $b= 1, \dots, c-1$, with the following 
properties:
\begin{itemize}
\item[(i)]
$\theta_b (\psi_1(C_b^{\nu_1}))$ is contained in a cube $\psi_0 (Q_b)$, where 
$Q_b$ is a cube of
$\eq \left( (1+2\nu_1)d_1 + \eps_1, A_{r_0}^{r_1} \right)$;
\item[(ii)]
among the cubes of $\eq \left( (1+2\nu_1)d_1 + \eps_1, A_{r_0}^{r_1} \right)$
different from $Q_1, \ldots, Q_{b-1}$,
the cube $Q_b$ is a cube that is closest to $B_{r_0}$;
\item[(iii)]
$\theta_b$ is supported in the domain
\[ 
\left( 
\widehat \cu_1 \setminus 
\psi_1 \left( \bigcup_{a=b+1}^{s_1} C_a^{\nu_1} \right)
\right)
\cup
\left( 
\cu_0'
\setminus     
\psi_0 \left( B_{r_0} \cup \bigcup_{a=1}^{b-1} Q_a \right) 
\right).
\]
\end{itemize}
In order to construct $\theta_c$, we first use the tree $\ct_1$ to construct 
a symplectomorphism $\tau_c$ of $\widetilde V_1$ whose support is contained in 
$\widehat V_1$ and is disjoint from $\bigcup_{a=c+1}^{s_1}C_a^{\nu_1}$ and which 
maps 
$C_c^{\nu_1}$ to the pilot cube $C_{\frak p}^{\nu_1}$.
This can be done as in the proof of Lemma~\ref{l:thetaj} by making use of 
Lemma~\ref{l:cover}\:(i) and the choice~\eqref{est:nud}.
This time, though, we may have to lift or lower blocking cubes by almost 
$\gd /2$, cf.~Figure~\ref{figure33}.
The smooth extension $\overline \tau_c$ of 
$\widetilde \psi_1 \circ \tau_c \circ \widetilde \psi_1^{-1}$ 
by the identity is supported in 
$\widehat \cu_1 \setminus \psi_1 \left( \bigcup_{a=c+1}^{s_1} 
C_a^{\nu_1}\right)$ and maps 
$\psi_1(C_c^{\nu_1})$ to $\widetilde \psi_1(C_{\frak p}^{\nu_1})$.
Since $C_{\frak p}^{\nu_1} \subset W_1 \cap \widetilde \psi_1^{-1} (\cu_0')$ 
and since $\cw_1 \subset \widetilde \cu_1$ is disjoint from $\psi_0 (B_{r_l})$,
we have that 
$\widetilde \psi_0^{-1} \circ \widetilde \psi_1 (C_{\frak p}^{\nu_1}) = C_{\frak 
p}^{\nu_1}$ is a cube in $B^{2n}(a_0') \setminus \overline{B_{r_l}}$.
After choosing $d_1$ yet smaller if necessary and in view of 
hypotheses~(i) and (ii) we therefore find a symplectomorphism $\gt_c$ 
supported in 
$U_0' \setminus \left( B_{r_0} \cup \bigcup_{b=1}^{c-1} Q_b \right)$ which maps 
$C_{\frak p}^{\nu_1}$ to a cube $Q_c$ of 
$\eq \left( (1+2\nu_1)d_1 + \eps_1, A_{r_0}^{r_1} \right)$ meeting (ii) with 
$b=c$.
The smooth extension $\overline \gt_c$ of 
$\psi_0 \circ \gt_c \circ \psi_0^{-1}$  
by the identity is supported in 
$\cu_0' \setminus \psi_0 \left( B_{r_0} \cup \bigcup_{b=1}^{c-1} Q_b \right)$ 
and maps 
$\psi_0 (C_{\frak p}^{\nu_1})$ to $\psi_0(Q_c)$.
The composition $\theta_c := \overline \gt_c \circ \overline \tau_c$ then meets 
properties~(i), (ii) and (iii).

Using these three properties of the maps 
$\theta_1, \dots, \theta_{s_1}$ as well as the inclusion 
$\Phi_0^j (\cs_0^j) \subset \psi_0 (B_{r_0})$ guaranteed by 
Proposition~\ref{p:Phi0} and the inclusion $\cs_1^j \subset \psi_1 (\es_1)$,
we see that
$\Phi_1^j := \theta_{s_1} \circ \dots \circ \theta_1$ 
is a symplectomorphism of $M$ whose support is disjoint from 
$\Phi_0^j (\cs_0^j) \cup \bigcup_{g=2}^l \cs_g^j$
and such that $\Phi_1^j (\cs_1^j) \subset \psi_0 (A_{r_0}^{r_1})$.

\s
Assume now by induction that for $h = 1, \dots, i-1$ 
we have already constructed symplectomorphisms  
$\Phi_h^j$ of $M$ whose support is disjoint from
\[
\bigcup_{g=0}^{h-1} \Phi_g^j \bigl( \cs_g^j \bigr) \cup
\bigcup_{g=h+1}^l \cs_g^j
\]
and such that
$\Phi_h^j \bigl( \cs_h^j \bigr) \subset \psi_0 \bigl( A_{r_{h-1}}^{r_h} \bigr)$. 
We are going to construct $\Phi_i^j$.
Recall from Step~3 that $\cs_i^j \subset \cu_i$.
As in the construction of $\Phi_1^j$ we consider a set of $s_i$ disjoint cubes
$\es_i = \bigcup C^{\nu_i} \subset V_i$ of width $(1+2 \nu_i)d_i$ containing the 
components of $\psi_i^{-1}(\cs_i^j)$.
The same reasoning and construction as for $\Phi_1^j$ shows that $d_i >0$ and 
$\nu_i >0$ (and for this $d_0, \dots, d_{i-1} >0$) can be chosen so small that 
for a suitable numbering of the cubes of $\es_i$ there are symplectomorphisms 
$\overline \tau_1, \dots, \overline \tau_{s_i}$ of $\widetilde \cu_i$ 
such that $\overline \tau_c$ is supported in 
\[
\widetilde \cu_i \setminus \left( \bigcup_{h=i+1}^l \cs_h^j \cup 
\bigcup_{a=c+1}^{s_i} \psi_i (C_a^{\nu_i})\right)
\]
and maps $\psi_i (C_c^{\nu_i})$ to a pilot cube 
$\widetilde \psi_i (C^{\nu_i}_{\frak p}) 
\subset \cw_i \cap \cu_{\underline i}'$.
Here, $\cw_i$ is the neighbourhood of $p_i \in \pp \cu_{\underline i}' \cap \pp 
\cu_i'$ on which $\widetilde \psi_{\underline i}^{-1} \circ \widetilde \psi_i$ 
is the identity.
In view of the estimate~\eqref{eq:rl} we can also assume that 
$\widetilde \cu_i$ is disjoint from $\psi_0 (B_{r_l})$, 
so that the supports of the $\overline \tau_c$ are also disjoint from $\psi_0 
(B_{r_l})$.

Recall now from Step~1 that all the sets $\cu_h'$ are non-empty and connected 
and have piecewise smooth boundary.
Moreover, $\cs_h^j \subset \cu_h$ for all $h$ and $\cu_g' \cap \cu_h = 
\emptyset$ if $g<h$ according to~\eqref{eq:uiuj}.
Therefore, 
\begin{equation}  \label{eq:disjoint}
\bigcup_{h=i}^l  \cs_h^j 
\quad \text{ is disjoint from} \quad
\bigcup_{g=0}^{i-1}  \cu_g'. 
\end{equation}
Let $0 < i_1 < \dots < \underline i < i$ be the branch from $\cu_0$ to $\cu_i$ 
in the rooted tree $\ct$ from Step~1.
Choosing the width $(1+2 \nu_i)d_i$ of $C^{\nu_i}_{\frak p}$ small enough,
we can use~\eqref{eq:disjoint} and the domains $\cu_{\underline i}', \dots, 
\cu_{i_1}'$ and the gates $\cw_i, \cw_{\underline i}, \dots, \cw_{i_2}$ and 
Lemma~\ref{l:trans} to construct a symplectomorphism $\theta$ of $M$
with support disjoint from 
$\bigl( \bigcup_{h=i}^l  \cs_h^j \bigr) \cup \psi_0 (B_{r_l})$,
and mapping $\widetilde \psi_i (C^{\nu_i}_{\frak p})$ to another pilot cube 
$\widetilde \psi_{i_1} (C^{\nu_i}_{\frak p}) \subset \cw_{i_1} \cap \cu_0'$.

Finally note that $\left| \es_i \right| < \bigl| A_{r_{i-1}}^{r_i} \bigr|$ and 
that for $d_i>0$ and $\nu_i>0$ small enough we find $\eps_i >0$ such that 
$\eq \bigl( (1+2\nu_i)d_i + \eps_i, A_{r_{i-1}}^{r_i} \bigr)$ 
contains at least $s_i$ cubes. 
As in the last step of the construction of $\Phi_1^j$ we therefore successively 
find $s_i$ symplectomorphisms $\overline \gt_c$ supported in 
$\cu_0' \setminus \psi_0 \left( B_{r_{i-1}} \cup \bigcup_{b=1}^{c-1}Q_b \right)$ 
and mapping $\widetilde \psi_{i_1}(C_{\frak p}^{\nu_i})$ to a cube $Q_c$ of 
$\eq \bigl( (1+2\nu_i)d_i + \eps_i, A_{r_{i-1}}^{r_i} \bigr)$.

After all, the symplectomorphism
\[
\Phi_i^j \,:\,=\, 
\left( \overline \gt_{s_i} \circ \theta \circ \overline \tau_{s_i} \right)
\circ \dots \circ
\left( \overline \gt_1 \circ \theta \circ \overline \tau_1 \right)
\]
has support disjoint from 
$\bigl( \bigcup_{h=i+1}^l  \cs_h^j \bigr) \cup 
\psi_0 \left( B_{r_{i-1}} \right)$
and maps $\psi_i (\es_i)$ into 
$\psi_0 \bigl( A^{r_i}_{r_{i-1}} \bigr)$.
Since 
$\bigcup_{g=0}^{i-1} \Phi_g^j(\cs_g^j) \subset \psi_0 (B_{r_{i-1}})$ by the 
induction hypothesis and since $\cs_i^j \subset \psi_i(\es_i)$,
the map $\Phi_i^j$ is as desired.
The proof of Proposition~\ref{p:Phih} is complete.
\proofend

In order to complete the proof of Theorem~\ref{t:ref}
we choose $d_0, \dots, d_l >0$ such that the conclusions of
Propositions~\ref{p:Phi0} and \ref{p:Phih} hold for each $j \in
\left\{ 1, \dots, k \right\}$, and we define the symplectomorphism
$\Phi^j$ of $M$ by
\[
\Phi^j \,=\, \Phi_h^j \circ \dots \circ \Phi_1^j \circ \Phi_0^j .
\]
In view of Propositions~\ref{p:Phi0} and \ref{p:Phih} and the
inclusion \eqref{inc:BAB} we then have 
\begin{eqnarray*}
\Phi^j \left( \cs^j \right) &=& \Phi^j \left( \bigcup_{h=0}^l \cs_h^j
                                \right) \\
           &=& \bigcup_{h=0}^l \Phi_h^j \left( \cs_h^j \right) \\
           &\subset& \psi_0 (B_{r_0}) \cup \bigcup_{h=1}^l \psi_0 \bigl( 
A_{r_{h-1}}^{r_h} \bigr) \\
           &\subset& \psi_0 \left( B^{2n}(a_0') \right) \\
           &\subset& \psi_0 \left( B^{2n}(a_0) \right) \\
           &=& \cb_0 .
\end{eqnarray*}
This and the identity~\eqref{id:MSj} imply that the $k$ Darboux charts
\[
\left( \Phi^j \right)^{-1} \circ \psi_0 \colon B^{2n}(a_0) \,\ra\, M
\]
cover $M$.
The proof of Theorem~\ref{t:ref} is finally complete, and so Theorem~1 is
also proved.
\proofend

\begin{remark}  \label{r:others}
{\rm
The method of the above proof can be used to obtain atlases with few charts in 
other situations.
For instance, one obtains the basic estimate $\B(M) \le \dim M + 1$ for closed 
connected manifolds proved in~\cite{L},
as well as the estimate $\C (M,\xi) \le \dim M + 1$ for the minimal number  
$C(M,\xi)$  of Darboux charts needed to cover a closed connected contact 
manifold $(M,\xi)$, see~\cite{NS}.
}
\end{remark}

\section{Variations of the theme}

\ni
Consider again a closed connected $2n$-dimensional symplectic manifold $(M, 
\go)$.
In the symplectic packing problem, one usually considers packings of
$(M, \go)$ by {\it equal}\:\! balls, see \cite{G1, MP, T2, B1, B2, Sch1, Sch2}.
In analogy to this, we define for each $a >0$ the invariant
\[
\SBa (M,\go) \,:\,=\, 
\min \left\{ k \mid M = \cb_1 \cup \dots \cup \cb_k \right\}
\]
where now each $\cb_i$ is the symplectic image $\gf_i \left( B^{2n}(a)\right)$ 
of {\it the same ball} $ B^{2n}(a)$, and 
where we set $\SBa (M,\go) = \infty$ if no such covering exists,
and we study the number
\[
\SBe (M,\go) \,:\,=\, \min_{a>0} \SBa (M,\go) .
\]
\begin{theorem}  \label{t:de} 
Let $(M, \go)$ be a closed connected symplectic manifold. Then
Theorem~1 holds with $\SB (M,\go)$ replaced by $\SBe (M,\go)$.
\end{theorem}
\proof
In the proof of Theorem~1 we have covered $(M, \go)$ by equal balls and
have thus proved Theorem~1 with $\SB (M,\go)$ replaced by $\SBe (M,\go)$.
\proofend

Clearly,
\begin{equation}  \label{e:sbsbe}
\SB (M,\go) \,\le\, \SBe (M,\go) .
\end{equation}
For every $a>0$ we denote by $\Emb \left( B(a), M \right)$ 
the space of symplectic embeddings of 
$\bigl( \overline{B^{2n}(a)}, \go_0 \bigr) \ha (M, \go)$ endowed with the
$C^\infty$-topology.
\begin{corollary}  \label{c:=}
Assume that $\gl (M, \go) \ge 2n+1$ or that 
$\Emb \left( B(a), M \right)$ is path-connected for all $a>0$.
Then $\SB (M,\go) = \SBe (M, \go)$. 
\end{corollary}

\proof
If $\gl (M, \go) \ge 2n+1$, then Theorem~1 and Theorem~\ref{t:de} 
yield $\SB (M, \go) = \gl (M, \go)$ and $\SBe (M, \go) = \gl (M, \go)$.

\s
Assume now that $\Emb \left( B(a), M \right)$ is path-connected for all $a>0$,
and choose $k = \SB (M, \go)$ symplectic embeddings
$\gf_i \colon \overline{B^{2n}(a_i)} \ha M$ such that 
$M = \bigcup_{i=1}^k \gf_i \left( B^{2n} (a_i) \right)$.
We choose $\eps >0$ so small that 
\[
M = \bigcup_{i=1}^k \gf_i \left( B^{2n} (a_i-\eps) \right) ,
\]
and set $a_i' = a_i-\eps$. 
We can assume that $a_1' = \max_i a_i'$. 
The identity $\SB (M,\go) = \SBe (M,\go)$ follows from
\begin{lemma} \label{l:=}
For each $i \ge 2$ there exists a symplectic embedding 
\[
\widetilde{\gf}_i \colon B^{2n} \left( a_1' \right) \ha M
\]
such that $\widetilde{\gf}_i |_{B^{2n}(a_i')} = \gf_i |_{B^{2n}(a_i')}$.
\end{lemma} 
\proof
By assumption, there exists a smooth family of symplectomorphisms
$\gf_i^t \colon B^{2n}(a_i) \ha M$ such that
\[
\gf_i^0 = \gf_1 |_{B^{2n}(a_i)}
\quad \text{ and } \quad
\gf_i^1 = \gf_i .
\]
Consider the subsets
\[
A = \bigcup_{t\in [0,1]} \{ t \} \times \gf_i^t \left( B^{2n}(a_i)
\right)
\quad \text{ and } \quad 
A' = \bigcup_{t\in [0,1]} \{ t \} \times \gf_i^t \left( B^{2n}(a'_i)
\right)
\]
of $[0,1] \times M$.
Since each set $\gf_i^t \left( B^{2n}(a_i) \right)$ is contractible,
there exists a smooth time-dependent Hamiltonian function
$ H \colon A\ra\, \RR$
generating the symplectic isotopy 
$\gf_i^t \circ \left( \gf_i^0 \right)^{-1} \colon 
\gf_1 \left( B^{2n}(a_i) \right) \ra M$.
By Whitney's Theorem there exists a smooth function 
$f \colon [0,1] \times M \ra [0,1]$ such that $f = 1$ on $A'$ and $f = 0$ on  $M 
\setminus A$. Let $\Phi \colon M \to M$ be the time-1-map of the flow generated 
by the Hamiltonian $fH$.
Then 
\[
\Phi  \,=\, \gf_i^1 \circ \left( \gf_i^0 \right)^{-1}
\quad \text{ on }\, \gf_1 \left( B^{2n}(a_i') \right) .
\]
For the embedding $\widetilde{\gf}_i \colon B^{2n}(a_i) \ha M$ defined by 
\[
\widetilde{\gf}_i \,:\,=\, \Phi \circ \gf_1 |_{B^{2n}(a_i)}
\]
we then find
\[
\widetilde{\gf}_i \,=\, 
\Phi \circ \gf_1  \,=\, 
\gf_i^1 \circ \left( \gf_i^0 \right)^{-1} \circ \gf_1 \,=\,
\gf_i^1 \circ \gf_1^{-1} \circ \gf_1 \,=\,
\gf_i^1 \quad \text{ on }\, B^{2n}(a_i') .
\]
The proof of Lemma~\ref{l:=} is complete, and so Corollary~\ref{c:=} is 
also proved.
\proofend

The spaces $\Emb \left( B(a), M \right)$ are known to be path-connected 
for all $a>0$ for $n=1$ and for a class of symplectic $4$-manifolds containing 
(blow-ups of) rational and ruled manifolds, see \cite{M}. 
No closed symplectic manifold is known 
for which $\Emb \left( B(a), M \right)$ is not path-connected for some $a>0$.
We thus ask
\begin{question}
{\rm
Is it true that $\SB (M, \go) = \SBe (M,\go)$ for every closed symplectic 
manifold $(M, \go)$?
}
\end{question}

\m
We next study the 
``symplectic Lusternik--Schnirelmann category'' $\S(M,\go)$ defined as
\[
\S (M, \go) \,=\, \min \left\{ k \mid M = \cu_1 \cup \dots \cup \cu_k \right\}
\]
where each $\cu_i$ is the image $\gf_i \left( U_i \right)$ of a
symplectic embedding $\gf_i \colon U_i \ra \cu_i \subset M$
of a bounded subset $U_i$ of $\left( \RR^{2n}, \go_0 \right)$ diffeomorphic to 
the open ball in $\RR^{2n}$.
\begin{theorem}  \label{t:S2n}
Let $(M, \go)$ be a closed connected $2n$-dimensional symplectic manifold. Then
$\S (M,\go) \le 2n+1$.
\end{theorem}
Theorem~\ref{t:S2n} will follow from a stronger result dealing with
covers by displaceable sets.
We say that a subset $\cu$ of $M$ is {\it displaceable}\,
if there exists an autonomous Hamiltonian function $H \colon M \ra \RR$
whose time-$1$-map $\gf_H$ displaces $\cu$, 
i.e., $\gf_H \left( \cu \right) \cap \cu = \emptyset$.
Define the invariant $\Sdis (M,\go)$ as
\[
\Sdis (M, \go) \,=\, \min \left\{ k \mid M = \cu_1 \cup \dots \cup
\cu_k \right\} 
\] 
where each $\cu_i$ is as in the definition of the invariant $\S (M,
\go)$ and is in addition displaceable.
Covers by such subsets $\cu_i$ play a role in the recent construction
of Calabi quasimorphisms on the group of Hamiltonian diffeomorphisms of
$(M, \go)$ in \cite{EP}, see also \cite{BEP}.
\begin{theorem}  \label{t:Sdis}
Let $(M, \go)$ be a closed $2n$-dimensional symplectic manifold. Then
$\Sdis (M,\go) \le 2n+1$.
\end{theorem}

Of course, $\B (M) \le \S (M, \go) \le \Sdis (M, \go)$.
Theorem~\ref{t:Sdis} thus implies Theorem~\ref{t:S2n}, 
and Proposition~1 and Theorem~\ref{t:Sdis} yield
\[
n+1 \,\le\, \cl (M)+1 \,\le\, \cat M \,\le\, \B (M) \,\le\, 
\S (M,\go) \,\le\, \Sdis (M,\go) \,\le\, 2n+1 
\]
and $\B (M) = \S (M,\go) = \Sdis (M,\go) = 2n+1$ if $[\go] |_{\pi_2(M)} =0$. 
For the $2$-sphere we have $2 = \S \left( S^2 \right) < 
\Sdis \left( S^2 \right) =3$.


\begin{question}  \label{q:BS}
{\rm
Is it true that $\B (M) = \S (M,\go)$ for every closed symplectic 
manifold $(M, \go)$?
}
\end{question}


\ni
{\it Proof of Theorem~\ref{t:Sdis}}:\,
Theorem~\ref{t:Sdis} is a consequence of the construction in the previous 
section and the following
\begin{proposition}  \label{p:eps}
For every $\eps >0$ there exists a symplectic embedding 
$\psi \colon (U, \go_0) \ha (M,\go)$ of a bounded subset $U$ of $\RR^{2n}$ 
diffeomorphic to a ball such that $\psi (U)$ is displaceable and
\[
\left| U \right| \,>\, \frac{\mu (M)}{2} -\eps .
\] 
\end{proposition}

Indeed, choose $\eps >0$ so small that 
\[
\frac{\mu (M)}{2} - \eps \,>\, \frac{\mu(M)}{2n+1} .
\]
For the set $\psi(U) \subset M$ guaranteed by Proposition~\ref{p:eps} we then
have
\[
\mu \left( \psi (U) \right) \,>\, \frac{\mu (M)}{2n+1} . 
\]  
Repeating the construction in the proof of Theorem~\ref{t:ref}  
with the ball $\cb = \gf \left( B^{2n}(a) \right)$ replaced by $\psi(U)$
and with $k=2n+1$, 
we find a cover $\{\cu_i\}$ of $M$ by $2n+1$ domains $\cu_i
\subset M$ which are diffeomorphic to balls and displaceable.

\begin{figure}[ht] 
 \begin{center}
  \psfrag{d}{$d$}
  \psfrag{3d}{$3d$}
  \psfrag{d2}{$\frac d2$}
  \psfrag{-d2}{$-\frac d2$}
  \psfrag{gd}{$\gd$}
  \psfrag{x1}{$x_1$}
  \psfrag{N}{$\cn$}
  \psfrag{U}{$U$}
 \leavevmode\epsfbox{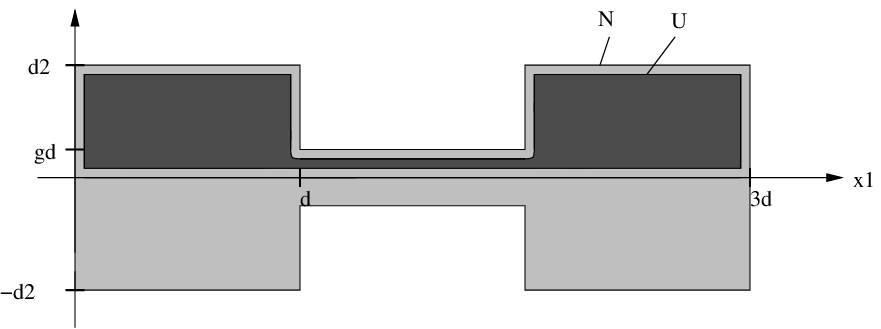}
\end{center}
\caption{The sets $\cn$ and $U$ for $k=1$.} 
 \label{figure34}
\end{figure}

\b
\ni
{\it Proof of Proposition~\ref{p:eps}}:\,
We fix $\eps >0$.
Let $k \in \NN$ and $d > \gd >0$.
For $j \in \NN \cup \{ 0 \}$ we denote by $\xi_{jd}$ the translation by
$jd$ in the $x_1$-direction and by $\eta_{-d/2}$ the translation by
$-d/2$ in the $y_1$-direction.
Consider the open subsets 
$C_j(d) = \xi_{2j} \left( \eta_{-d/2} \bigl( \left] 0,d \right[^{2n}
\bigr) \right)$ 
and 
\[
\cn \left( k,d,\gd \right) \,=\, \coprod_{j=0}^k C_j(d) \cup 
\left( \left] 0,(2k+1)d\right[ \times \left] -\gd,\gd \right[^{2n-1} \right)
\]
of $( \RR^{2n}, \omega_0 )$.
Figure~\ref{figure34} illustrates a set $\cn \left( k,d,\gd \right) \subset 
\RR^{2n}$
for $k=1$.

%
%

\ni
According to \cite[Section~6.1]{Sch1} there exist $k$, $d$ and $\gd$ and
a symplectic embedding $\psi \colon \cn \left( k,d,\gd \right) \ha
(M,\go)$ such that 
\begin{equation}  \label{e:kme}
\left| \coprod_{j=0}^k C_j(d) \right|
\,>\, \mu (M) - \eps .
\end{equation}
Set $\cn^+ ( k,d,\gd ) = \cn \left( k,d,\gd \right) 
\cap \left\{ y_1 > 0 \right\}$, and denote by $\pp \cn^+ ( k,d,\gd )$
the boundary of this set.
For $\nu>0$ we set 
\[
U_\nu \,=\, \left\{ z \in \cn^+ ( k,d,\gd ) \mid 
              \dist \left( z, \pp \cn^+ (k,d,\gd) \right) > \nu
              \right\},
\] 
cf.~Figure~\ref{figure34}.
For $\nu < \gd /2$ the set $U_\nu$ is connected and diffeomorphic to a
ball.
In view of~\eqref{e:kme} we can choose $\nu < \gd/2$ so small that
\[
\left| U_\nu \right| \,>\, \frac{\mu (M)}{2} - \eps .
\] 
For such a choice of $k$, $d$, $\gd$ and $\nu$ we abbreviate 
$\cn = \cn \left( k,d,\gd \right)$ and $U = U_\nu$.
We shall construct a Hamiltonian isotopy $\gf_t$ 
of $\RR^{2n}$ which is generated by an autonomous Hamiltonian function
with support in $\cn$ and such that $\gf_1 (U) \cap U = \emptyset$.
The autonomous Hamiltonian diffeomorphism $\Phi$ of $(M,\go)$ defined by
\[
\Phi (z) \,=\, 
         \left\{
           \begin{array}{ll}   
                \psi \circ \gf_1 \circ \psi^{-1} (z) & 
                     \text{ if } z\in \psi \left( \cn \right) \\[.2em]
                 z & \text{ if } z\notin \psi \left( \cn \right) 
           \end{array}   \right. 
\]
then displaces $\psi \left( U \right)$.
Note that the images of $\cn$ and $U$ under the projection
$\pi \colon \RR^{2n} \ra \RR^2$, 
$(x_1,y_1, \dots, x_n,y_n) \mapsto (x_1,y_1)$,
look as in Figure~\ref{figure34}.
In order to construct the Hamiltonian isotopy $\gf_t$, 
we choose a smooth function $f \colon \RR \ra \RR$ such that on 
$\left] 0, (2k+1)d \right[$ the graph of $f$ is 
contained in $\pi \left( \cn \right)$ and lies above $\pi (U)$.
Then the Hamiltonian function $H \colon \RR^{2n} \ra \RR$ defined by
\[
H (x_1,y_1,x_2, \dots, y_n) \,=\, - \int_0^{x_1} f(s)\,ds 
\]
generates the isotopy 
\[
\phi_t \colon (x_1,y_1,x_2, \dots, y_n) \,\mapsto\, 
\left (x_1,y_1 - t f(x_1), x_2, \dots, y_n \right), \quad t \in [0,1],
\]
which satisfies $\phi_t (U) \subset \cn$ for all $t \in [0,1]$ and
$\phi_1(U) \cap U = \emptyset$.
Choose now a smooth function $h \colon \RR^{2n} \ra [0,1]$ which is
equal to $1$ on 
$\bigcup_{t \in [0,1]} \phi_t (U)$  
and vanishes outside $\cn$.
The Hamiltonian isotopy $\gf_t$ generated by the Hamiltonian function
$hH$ is then as required.
\proofend

\section{Proof of Proposition~1}  \label{s:p1}

\ni
Since $(M, \go)$ is symplectic, $[\go]^n \neq 0$, and so
$n+1 \le \cl (M) +1$. The first statement in Proposition~1
follows from this estimate and from~\eqref{e:F}.

A main ingredient in the remainder of the proof is the following theorem
of W.~Singhof, who thoroughly studied the relation between $\B (M)$ and $\cat 
M$. 
Recall that a topological space $X$ is said to be $p$-connected if it is 
path-connected and its homotopy groups $\pi_i (X)$ vanish for $1 \leq i \le p$.
\begin{theorem} \label{t:singhof}
{\rm (Singhof, \cite[Corollary (6.4)]{S})}
Let $M^m$ be a closed smooth $p$-connected manifold with $m \ge 4$ and
$\cat M \ge 3$. Then
\begin{itemize}
\item[(a)]
$\B (M) \, = \, \cat M 
\quad \text{if } \cat M \, \ge \, \dfrac{m+p+4}{2(p+1)}$;
\item[(b)]
$\B (M) \, \le \, \left\lceil \dfrac{m+p+4}{2(p+1)} \right\rceil 
\quad \text{if } \cat M \, < \, \dfrac{m+p+4}{2(p+1)}$ .
\end{itemize}
$($Here, $\lceil x \rceil$ denotes the minimal integer which is
greater than or equal to $x$.$)$
\end{theorem} 

Since we consider only symplectic manifolds, 
the assumptions $\dim M \ge 4$ and $\cat M \ge 3$ in 
Theorem~\ref{t:singhof} can be dropped. 
Indeed, if $\dim M = 2$, it is easy to
see that we are in the situation of (a) in Theorem~\ref{t:singhof}; and
if $\cat M =2$, then $\frac 12 \dim M \le \cl (M) +1 \le \cat M =2$ yields
$\dim M =2$.

\s
(i)
If $M$ is simply connected, then $\cat M \le n+1$, see \cite{J1}, 
and so $\cat M = n+1$.
This and again $p \ge 1$ show that we are in the situation of 
Theorem~\ref{t:singhof}\,(a), so $\B(M) = \cat M$.

(ii) 
It has been proved in \cite{RO} 
that $[\go]|_{\pi_2(M)} =0$ implies $\cat M = 2n+1$, and so the claim 
follows together with
$\B(M) \le 2n+1$.

(iii) 
As we remarked above, $\B(M) = \cat M$ if $n=1$. 
So let $n \ge 2$ and assume that $\B(M) > \cat M$. By (i) we have $p=0$. The 
claim
now readily follows from Theorem~\ref{t:singhof}.
\proofend

\begin{remarks}\
{\rm
{\bf 1.} 
The inequality $\cl (M) +1 \le \cat M$ can be strict:
For the Thurston--Kodaira manifold described in \cite[Example~3.8]{MS}
we have $\pi_2(M) =0$ and hence $\cat M =5$, but $\cl (M) =3$, see \cite{R}.
More generally, $\cl(M)+1 < \cat M = \dim M +1$ for any symplectic 
non-toral nilmanifold, see \cite{RT}. 

\s
{\bf 2.}
It follows from \cite[Prop.~13]{LSV} and \cite[Prop.~3.6]{CLOT}
that there do exist closed smooth manifolds with $\cat M < \B (M)$.
No symplectic examples are known, however.
}
\end{remarks}

\begin{examples}\
{\rm
{\bf 1.}
Examples of closed symplectic manifolds $(M^{2n},\go)$ which are simply 
connected and hence have $\B(M)=n+1$ are hyperplane sections of 
$(\CP^{n+1},\go_{SF})$ with $n \ge 2$. 
Many more such examples can be found in~\cite{Go}.

\s
{\bf 2.}
If $\left( M^{2n}, \go \right)$ admits a Riemannian metric with nonnegative 
Ricci curvature and has infinite fundamental group, then
\[
\cat M \,\ge\, n+1 + \frac{b_1(M)}{2}
\quad \text{ and } \quad b_1(M) >0, 
\]
see \cite[Theorem~4.3]{O}.
In particular, $\cat M \ge n+2$, and so $\cat M = \B (M)$ by 
Proposition~1\:(iii).

\s
{\bf 3.}
Assume that the homomorphism
$[\go]^{n-1} \colon H^1 (M;\RR) \ra H^{2n-1} (M;\RR)$ 
(multiplication by the class $[\go]^{n-1}$) is a non-zero map. 
K\"ahler manifolds with $H^1(M;\RR) \neq 0$ have this property.
Using Poincar\'e duality we see that
$\cl (M) \ge n+1$, and so $n+2 \le \cat M = \B(M)$.
}
\end{examples}

\section{Examples}  \label{s:ex}

\ni
In this section we compute or estimate the number $\SB (M, \go)$ for
various closed symplectic manifolds $(M, \go)$. 
In view of the estimate
\[
\SB (M^{2n}, \go) \,\ge\, 
\Gamma (M, \go) \,=\, 
\left\lfloor \frac{\Vol (M, \go)}{\frac{1}{n!} \left( \Gr(M,\go) \right)^n} 
\right\rfloor +1 
\]
from Theorem~1 
and in view of Proposition~1, understanding $\SB (M, \go)$ is
often equivalent to understanding the Gromov width $\Gr (M, \go)$. 
Our list of examples therefore resembles the list of closed symplectic manifolds 
whose Gromov width is known,
\cite{B1,B2,BC,Ji,KT,LM1,LM3,Lu,M-variants,MP,MSl,Sch1,Sib}.

An important tool for obtaining upper bounds of the Gromov width in many 
examples is Gromov's Nonsqueezing Theorem.
The proof of the following general version makes use of the existence of 
Gromov--Witten invariants for arbitrary closed symplectic manifolds,
see \cite[Section~9.3]{MS2} and \cite[Proposition~1.18]{M-variants}.

\begin{NStheorem} \label{t:NS}
For any closed symplectic manifold $(M,\go_M)$,
\[
\Gr \left( M \times S^2, \go_M \oplus \go_{S^2} \right) \, \le \, 
\int_{S^2}\go_{S^2} .
\]
\end{NStheorem} 

\ni
For generalizations of this result we refer to \cite[Section~9.3]{MS2},
\cite[Remark~9.3.7]{Sch1} and \cite{Lu}.

\m
We shall also frequently use the following well-known fact.
\begin{lemma} \label{l:GS}
{\rm (Greene--Shiohama, \cite{GS})}
Let $U$ and $V$ be bounded domains in $\left( \RR^2, dx \wedge dy
\right)$ which are diffeomorphic and have equal area. 
Then $U$ and $V$ are symplectomorphic.
\end{lemma} 

\m
\ni
{\bf 1. Surfaces.}
A closed $2$-dimensional symplectic manifold is a closed oriented surface
equipped with an area form.

\begin{corollary}  \label{c:surface}
Let $(\Sigma_g, \gs)$ be a closed oriented surface of genus $g$ 
with area form $\gs$. Then
\[
\SB (\Sigma_g, \gs) =                       
    \left\{
           \begin{array}{l}   
                  2 \;\, \text{ if } g = 0,  \\ [0.2em]
                  3 \;\, \text{ if } g \ge 1 .
           \end{array}   \right. 
\]  
\end{corollary}

\proof
In view of Lemma~\ref{l:GS} we have $\SB (\Sigma_g, \gs) = \B(\Sigma_g)$, and so 
the
corollary follows in view of Proposition~1.
\proofend

\m
\ni
{\bf 2. Minimal ruled 4-manifolds.}
As before we denote by $\Sigma_g$ the closed oriented surface of genus $g$.
There are exactly two orientable $S^2$-bundles with base 
$\Sigma_g$, namely the trivial bundle 
$\Sigma_g \times S^2 \ra \Sigma_g$ and the
nontrivial bundle $\Sigma_g \ltimes S^2 \ra \Sigma_g$, 
see \cite[Lemma~6.9]{MS}.

\m
\ni
{\it a) Trivial $S^2$-bundles.}
Fix area forms $\gs_{\Sigma_g}$ and $\gs_{S^2}$ of area 1 on $\Sigma_g$
and $S^2$, respectively.
By the work of Lalonde--Mc\,Duff and Li--Liu
every symplectic form on $\Sigma_g \times S^2$ is diffeomorphic to 
$a \gs_{\Sigma_g} \oplus b \gs_{S^2}$ for some $a,b >0$ (see \cite{LM2}).
We abbreviate $\Sigma_g(a) = (\Sigma_g, a \gs_{\Sigma_g})$ and
$S^2(b) = (S^2, b \gs_{S^2})$.

\begin{figure}[ht] 
 \begin{center}
  \psfrag{1}{$1$}
  \psfrag{32}{$\frac{3}{2}$}
  \psfrag{2}{$2$}
  \psfrag{3}{$3$}
  \psfrag{4}{$4$}
  \psfrag{5}{$5$}
  \psfrag{6}{$6$}
  \psfrag{7}{$7$}
  \psfrag{8}{$8$}
  \psfrag{b}{$\tfrac{a}{b}$}
  \psfrag{SB}{$\SB \left( S^2(a) \times S^2(b) \right)$}
  \leavevmode\epsfbox{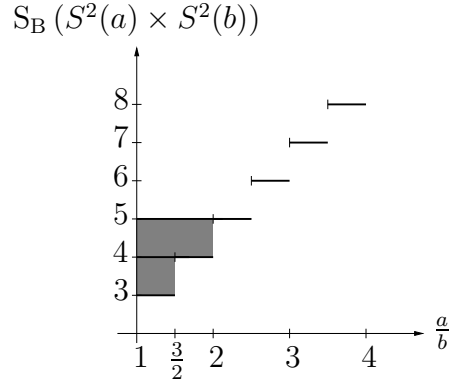}
 \end{center}
 \caption{What is known about $\SB \left( S^2(a) \times S^2(b) \right)$.} 
 \label{figure11}
\end{figure}

\begin{figure}[ht] 
 \begin{center}
  \psfrag{0}{$0$}
  \psfrag{1}{$1$}
  \psfrag{2}{$2$}
  \psfrag{3}{$3$}
  \psfrag{4}{$4$}
  \psfrag{5}{$5$}
  \psfrag{6}{$6$}
  \psfrag{7}{$7$}
  \psfrag{8}{$8$}
  \psfrag{b}{$\tfrac{a}{b}$}
  \psfrag{SB}{$\SB \left( \Sigma_g(a) \times S^2(b) \right)$}
  \leavevmode\epsfbox{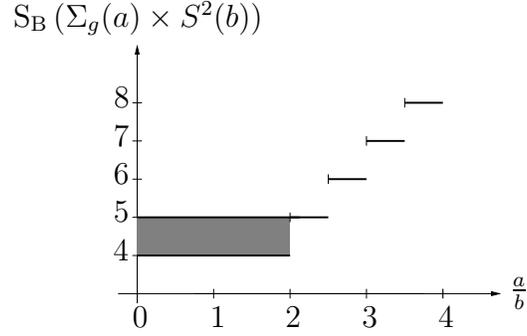}
 \end{center}
 \caption{What is known about $\SB \left( \Sigma_g(a) \times S^2(b)
 \right)$.} 
 \label{figure12}
\end{figure}

\begin{corollary}  \label{c:trivial}
For $S^2(a) \times S^2(b)$ with $a \ge b 
>0$
we have 
\[
\SB \left( S^2(a) \times S^2(b) \right)                      
    \left\{
      \begin{array}{lll}   
           \in \{ 3,4,5\} & \text{ if}&  1 \le \tfrac{a}{b} <
                                               \tfrac{3}{2},  \\[0.2em]
           \in \{ 4,5 \} & \text{ if}&  \tfrac{3}{2} \le
                                               \tfrac{a}{b} < 2,   \\[0.2em]
           = \left\lfloor \tfrac{2a}{b} \right\rfloor +1 & \text{ if}&
                                         \tfrac{a}{b} \ge 2, 
      \end{array}   \right. 
\]  
and for $\Sigma_g(a) \times S^2(b)$ with 
$g \ge 1$ and $a,b >0$ we have
\[
\SB \left( \Sigma_g(a) \times S^2(b) \right)                       
    \left\{
      \begin{array}{lll}   
             \in \{ 4,5 \} & \text{ if}&  0 < \tfrac{a}{b} < 2,   \\[0.2em]
             = \left\lfloor \tfrac{2a}{b} \right\rfloor +1 & \text{ if}&
                                              \tfrac{a}{b} \ge 2 .
      \end{array}   \right. 
\]  
\end{corollary}
\ni
The result is illustrated in Figures~\ref{figure11} and \ref{figure12}.

%
%

%
%
%

\b
\proof
Proposition~1\:(i) yields $\B(S^2 \times S^2) =3$. 
Moreover, Gromov's Nonsqueezing Theorem~\ref{t:NS} implies that        
$\Gr (S^2(a) \times S^2(b)) = b$, and so 
\[
\Gamma \left( S^2(a) \times S^2(b) \right) = 
\left\lfloor \tfrac{2a}{b} \right\rfloor +1.
\]
The first half of the corollary now follows from Theorem~1.

For any two path-connected $CW$-spaces $X$ and $Y$ it holds that
\[
\cat ( X \times Y) \,<\, \cat X + \cat Y, 
\]
see \cite{J1}.
This and $\cl (\Sigma_g \times S^2)+1=4$ show that
$\cat (\Sigma_g \times S^2) = 4$, and so $\B (\Sigma_g \times S^2) = 4$
in view of Proposition~1\:(iii).
Moreover, it follows from Theorem~6.1.A in \cite{B1} that
\[
\Gamma  \left( \Sigma_g(a) \times S^2(b) \right) \,=\, 
\left\lfloor \max \left\{ 1, \tfrac{2a}{b} \right\} \right\rfloor +1 .
\]
The second half of Corollary~\ref{c:trivial} now follows from Theorem~1.
\proofend

\ni
{\it b) Nontrivial $S^2$-bundles.}
Let $A \in H_2(\Sigma_g \ltimes S^2; \ZZ)$ be the class of a section
with self intersection number $-1$, and let $F$ be the
homology class of the fiber. We set $B= A+ \tfrac{1}{2} F$. 
Then $\{ F, B \}$ is a basis of $H_2(\Sigma_g \ltimes S^2; \RR)$.
For $a,b >0$ we fix a representative $\go_{ab}$ of the Poincar\'e dual of 
$aF+bB$.
By \cite[Theorem~6.11]{MS} and the work of Lalonde--Mc\,Duff and Li--Liu 
(see \cite{LM2}),

\m
1. Every symplectic form on $S^2 \ltimes S^2$ is diffeomorphic to 
$\omega_{ab}$ for some $a > \tfrac{b}{2} > 0$.

2. Every symplectic form on $\Sigma_g \ltimes S^2$, $g \ge 1$, 
is diffeomorphic to 
$\omega_{ab}$ for some $a,b > 0$.

\begin{corollary}  \label{c:nontrivial}
For $(S^2 \ltimes S^2, \go_{ab})$ with $a> \tfrac{b}{2} >0$
we have 
\[
\SB \left( S^2 \ltimes S^2, \go_{ab} \right)                      
    \left\{
      \begin{array}{lll}   
             \in \{ 3,4,5 \} & \text{ if}& \tfrac{1}{2} < \tfrac{a}{b} <
                                                   \tfrac{3}{2},   \\[0.2em]
             \in \{ 4,5 \} & \text{ if}&  \tfrac{3}{2} \le
                                                  \tfrac{a}{b} < 2,   \\[0.2em]
             = \left\lfloor \tfrac{2a}{b} \right\rfloor +1 & \text{ if}&
                                              \tfrac{a}{b} \ge 2 , 
      \end{array}   \right. 
\]  
and for 
$\left( \Sigma_g \ltimes S^2, \go_{ab} \right)$ with $g \ge 1$ and $a,b>0$ we 
have
\[
\SB \left( \Sigma_g \ltimes S^2, \go_{ab} \right)                       
    \left\{
      \begin{array}{lll}   
             \in \{ 4,5 \} & \text{ if}&  0 < \tfrac{a}{b} < 2,   \\[0.2em]
             = \left\lfloor \tfrac{2a}{b} \right\rfloor +1 & \text{ if}&
                                              \tfrac{a}{b} \ge 2 .
      \end{array}   \right. 
\]  
\end{corollary}

\begin{figure}[ht] 
 \begin{center}
  \psfrag{12}{$\frac{1}{2}$}
  \psfrag{1}{$1$}
  \psfrag{32}{$\frac{3}{2}$}
  \psfrag{2}{$2$}
  \psfrag{3}{$3$}
  \psfrag{4}{$4$}
  \psfrag{5}{$5$}
  \psfrag{6}{$6$}
  \psfrag{7}{$7$}
  \psfrag{8}{$8$}
  \psfrag{b}{$\tfrac{a}{b}$}
  \psfrag{SB}{$\SB \left( S^2(a) \ltimes S^2(b) \right)$}
  \leavevmode\epsfbox{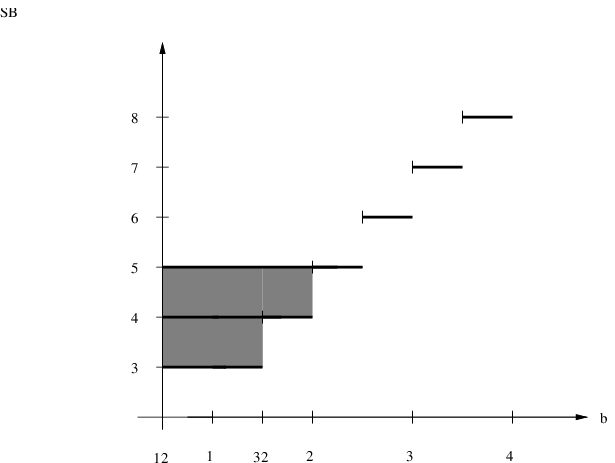}
 \end{center}
 \caption{What is known about $\SB \left( S^2 \ltimes S^2, \go_{ab} \right)$.} 
 \label{figure11b}
\end{figure}

\begin{figure}[ht] 
 \begin{center}
  \psfrag{0}{$0$}
  \psfrag{1}{$1$}
  \psfrag{2}{$2$}
  \psfrag{3}{$3$}
  \psfrag{4}{$4$}
  \psfrag{5}{$5$}
  \psfrag{6}{$6$}
  \psfrag{7}{$7$}
  \psfrag{8}{$8$}
  \psfrag{b}{$\tfrac{a}{b}$}
  \psfrag{SB}{$\SB \left( \Sigma_g \ltimes S^2, \go_{ab} \right)$} 
  \leavevmode\epsfbox{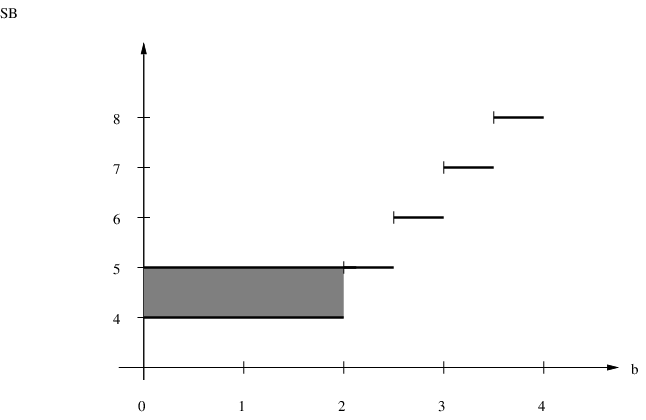}
 \end{center}
 \caption{What is known about $\SB \left( \Sigma_g \ltimes S^2, \go_{ab} 
\right)$.} 
 \label{figure12b}
\end{figure}

\ni
The result is illustrated in Figures~\ref{figure11b} and \ref{figure12b}.

%
%

%
%
%

\b
\proof
Since $S^2 \ltimes S^2$ is simply connected, 
$\B(S^2 \ltimes S^2) =3$ in view of Proposition~1\:(i). 
Moreover, Biran's work~\cite{B1} implies
\[
\Gamma\left (S^2 \ltimes S^2, \go_{ab} \right) = \left\lfloor \tfrac{2a}{b} 
\right\rfloor +1 ,
\]
see \cite{Sch2}.
The first half of the corollary now follows from Theorem~1.

Using the Leray--Hirsch Theorem,
we find that $\cl \left( \Sigma_g \ltimes S^2 \right) =3$, and so
$\cat \left( \Sigma_g \ltimes S^2 \right) \ge 4$.
On the other hand, $\Sigma_g \ltimes S^2$ having a section, it is not
hard to see that
$\cat \left( \Sigma_g \ltimes S^2 \right) \le 4$ (cf.\ the proof of
Proposition~3.3 in \cite{S1}).
In view of Proposition~1\:(iii) we conclude that
$\B (\Sigma_g \ltimes S^2) = 4$.
Moreover, it has been computed in \cite{Sch2} that
\[
\Gamma  \left( \Sigma_g \ltimes S^2, \go_{ab} \right) \,=\, 
\left\lfloor \max \left\{ 1, \tfrac{2a}{b} \right\} \right\rfloor +1 .
\]
The second half of the corollary now follows from Theorem~1.
\proofend

\b
\ni
{\bf 3. Products of higher genus surfaces.}
As before we denote by $\Sigma_g$ the closed oriented surface of genus
$g$. 
In view of the previous example we assume $g \ge 1$.
By a theorem of Moser \cite{Mo}, any two area forms on $\Sigma_g$ of 
total area $a$ are diffeomorphic. We write $\Sigma_g(a)$ for this
symplectic manifold.

\begin{corollary}\  \label{c:product}
\begin{itemize}
\vspace{0.2em}
\item[(i)] 
$\SB \left( \Sigma_1(a) \times \Sigma_g(b) \right) \;\!\!= 5$\; if $\tfrac{a}{b} 
< 
\tfrac{5}{2}$.
\vspace{0.2em}
\item[(ii)]
$\SB \left( \Sigma_g (a) \times \Sigma_h(b) \right) = 5$\; if
$\tfrac{2}{5} < \tfrac{a}{b} < 
\tfrac{5}{2}$.
\end{itemize}
\vspace{0.2em}
\end{corollary}

\proof
By Proposition~1\:(ii) we have that
\[
\B \left( \Sigma_g \times \Sigma_h \right) \,=\, 5 \quad \text{for all } g,h
\ge 1.
\]
Using Lemma~\ref{l:GS}
we see that the discs $B^2(a)$ and $B^2(b)$ symplectically embed into
$\Sigma_g(a)$ and $\Sigma_h(b)$, respectively. 
Therefore, the ball $B^4( \min (a,b) ) \subset B^2(a) \times B^2(b)$
symplectically embeds into $\Sigma_g(a) \times \Sigma_h(b)$, 
and so 
\[
\Gamma \bigl( \Sigma_g(a) \times \Sigma_h(b) \bigr) \le 5
\quad \text{whenever } \tfrac{2}{5} < \tfrac{a}{b} < \tfrac{5}{2} .
\]
Claim (ii) now follows from Theorem~1.

We prove Claim~(i) following~\cite{Ji}. For each $c>0$ we
consider the rectangle 
\[
R(c) = \left\{ (x,y) \in \RR^2 \mid 0<x<1, \; 0<y<c \right\},
\]
and the linear symplectic map
\begin{eqnarray*}
\gf \colon  (R(c) \times R(c), \go_0)  &  \ra   &
                (\RR^2 \times \RR^2, \go_0)   \\
      (x_1,y_1,x_2,y_2) & \mapsto & (x_1+y_2,y_1, -y_2, y_1+x_2)
\end{eqnarray*}
where $\go_0 = dx_1 \wedge dy_1 + dx_2 \wedge dy_2$.
Let $T^2 = \left( \RR^2 / \ZZ^2, dx_1 \wedge dy_1 \right)$ be the
standard symplectic torus.
Then the projection $p \colon \left( \RR^2, dx_1 \wedge dy_1 \right) \ra
T^2$ is symplectic, and so the composition
\[
(p \times id) \circ \gf \colon R(c) \times R(c) \ra T^2 \times \RR^2
\]
is also symplectic.
It is easy to see that this map is an embedding and
that
%
\[
\bigl( (p \times id) \circ \gf \bigr) \bigl( R(c) \times R(c) \bigr) 
\,\subset\, T^2 \:\! \times \;\! ]-c,0[ \times ]0, c+1[. 
\]
In view of Lemma~\ref{l:GS} the ball $B^4(c)$ symplectically embeds into
$R(c) \times R(c)$, and $]-c,0[ \times ]0, c+1[$ symplectically embeds into 
$\Sigma_g (c(c+1))$. 
We conclude that the ball $B^4(c)$ symplectically embeds
into $\Sigma_1(1) \times \Sigma_g (c(c+1))$ for each $c>0$, i.e., 
\[
\Gr \left( \Sigma_1(1) \times \Sigma_g(d) \right) 
\,\ge\, \tfrac{1}{2} \left( \sqrt{4d+1} -1 \right) \quad \text{for each } d>0 .
\]
This estimate and a computation yield
\[
\Gamma \left( \Sigma_1(a) \times \Sigma_g(b) \right) \,=\,
\Gamma \left( \Sigma_1(1) \times \Sigma_g \left( \tfrac{b}{a} \right) \right)
\,\le\ 5 \quad \text{whenever } \tfrac{a}{b} < \tfrac{9}{10}.
\]
Now, the already proved Claim~(ii) and Theorem~1 imply Claim~(i).
\proofend

\begin{remarks}\ 
{\rm
{\bf 1.}
Assume that $g \ge 1$, $h \ge 2$ and $\tfrac{a}{b} \ge \tfrac{5}{2}$.
The method used in the proof of (ii) in Corollary~\ref{c:product} only
yields the linear estimate
\[
\SB \left(  \Sigma_g (a) \times \Sigma_h (b) \right)
\,\le\, \left\lfloor \tfrac{2a}{b} \right\rfloor +1 .
\]
A variant of the method used in the proof of (i), however, yields the
estimate 
\[
\SB \left(  \Sigma_g (a) \times \Sigma_h (b) \right)
\,\le\, C(h) \, \frac{\tfrac{a}{b}}{\left( \log \tfrac{a}{b} \right)^2}
\] 
where $C(h)>0$ is a constant depending only on $h$ (see \cite{Ji}) .

\s
{\bf 2.}
Symplectic structures on torus bundles over closed orientable surfaces were 
studied in \cite{Ge,Kahn,Th,Wa},
but their Gromov widths are not known.
}
\end{remarks}

\b
\ni
{\bf 4. Complex projective space.}
Let $\CC \PP^n$ be complex projective space and let
$\go_{SF}$ be the unique $\mbox{U}(n+1)$-invariant K\"ahler form on 
$\CC \PP^n$ whose integral over $\CC \PP^1$ equals $1$. 

\begin{corollary}  \label{c:cpn}
$\SB \left( \CC \PP^n, \go_{SF} \right) \,=\, n+1$.
\end{corollary}

\proof
In view of Proposition~1\:(ii) we have 
\[
\SB \left( \CC \PP^n, \go_{SF} \right) \,\ge\, \B \left( \CC \PP^n
\right) \,\ge\, n+1 .
\]
On the other hand, we define for
$0 \le i \le n$ maps $f_i \colon B^{2n}(1) \ra \CC \PP^n$
by 
\begin{equation} \label{e:chartCPn}
f_i \colon \mathbf{z} = (z_1, \dots, z_n) \, \mapsto \, \left[ z_1 : \ldots :
z_{i-1} : \sqrt{1 - |\mathbf{z}|^2} : z_{i+1} : \ldots : z_n \right].
\end{equation}
It is well known that $f_i$ is a symplectomorphism between $B^{2n}(1)$
and $\CC \PP^n \setminus S_i$, where $S_i = \{ [ u_1 : \ldots : u_{i-1} :
0 : u_{i+1} : \ldots : u_n ] \} \cong \CC \PP^{n-1}$ is the $i$-th
coordinate hypersurface (see e.g.\ \cite{K}). 
Since 
\[
\CC \PP^n \,\subset\, \bigcup_{i=0}^n f_i \left( B^{2n}(1) \right),
\] 
we conclude that also 
$\SB (\CC \PP^n, \go_{SF}) \le n+1$, and so the corollary follows.
\proofend

\begin{remark} \label{r:taubes}
{\rm
By a theorem of Taubes, \cite{T}, any symplectic form on $\CP^2$ is
diffeomorphic to $a \;\!\go_{SF}$ for some $a>0$.
In view of Corollary~\ref{c:cpn} we thus have
\[
\SB( \CP^2, \go) = 3 \quad \text{\it for any symplectic form } \go
\text{ \it on } \CC 
P^2.
\]
}
\end{remark}

\b
\ni
{\bf 5. Complex Grassmann manifolds.}
Let $G_{k,n}$ be the Grassmann manifold of $k$-planes in $\CC^n$, and
let
$\gs_{k,n}$ be the standard K\"ahler form on $G_{k,n}$ normalized such
that $\gs_{k,n}$ is Poincar\'e dual to the generator of $H_2(G_{k,n};
\ZZ) = \ZZ$.
Since $(G_{n-k,n}, \gs_{n-k,n}) = (G_{k,n}, \gs_{k,n})$, we can assume
that
\[
k \in \left\{ 1, \dots, \left\lfloor \tfrac{n}{2} \right\rfloor \right\} .
\]
We define the number $p_{k,n}$ by
\begin{equation}  \label{d:pkn}
p_{k,n} \,=\, 
\frac{ (k-1)! \cdots 2! \, 1! \cdot (k(n-k))!}{(n-1)! \cdots (n-k+1)! \,
(n-k)!} .
\end{equation} 
Notice that $p_{k,n} = \deg (p(G_{k,n}))$ where
\[
p \colon G_{k,n} \ha \CP^{\binom n k -1}
\]
is the Pl\"ucker map \cite[Example 14.7.11]{Fu}, and so $p_{k,n}$ is
indeed an integer.
Since $(G_{1,n}, \gs_{1,n}) = ( \CP^{n-1}, \go_{SF})$, we assume $k
\ge 2$.

\begin{corollary}\  \label{c:grassmann}
\begin{itemize}
\item[1)]
$\SB (G_{2,4}, \gs_{2,4}) \in \{ 5, 6 \}$, \\
$\SB (G_{2,5}, \gs_{2,5}) \in \{ 7, 8, 9, 10 \}$, \\
$\SB (G_{2,n}, \gs_{2,n}) = p_{2,n}+1$ for all $n \ge 6$,
\vspace{0.2em}
\item[2)]
$\SB (G_{k,n}, \gs_{k,n}) = p_{k,n}+1$ for all $k \ge 3$.
\end{itemize}
\end{corollary}

\proof
Since $G_{k,n}$ is simply connected and since
\begin{equation}  \label{idg:d}
\dim G_{k,n} \,=\, 2k(n-k),
\end{equation}
we read off from Proposition~1\:(i) that
\begin{equation}  \label{idg:B}
\B \left( G_{k,n} \right) \,=\, k(n-k) +1 .
\end{equation}
Moreover, 
\begin{equation}  \label{id:vol}
\Vol \left( G_{k,n}, \gs_{k,n} \right) \,=\, 
\frac{p_{k,n}}{(k (n-k))!}
\end{equation}
(see \cite[Example 14.7.11]{Fu}), and it has been proved in \cite{KT,Lu} that
\[
\Gr \left( G_{k,n}, \gs_{k,n} \right) \,=\, 1 .
\]
Therefore,
\begin{equation}  \label{idg:G}
\Gamma \left( G_{k,n}, \gs_{k,n} \right) \,=\, p_{k,n} +1 .
\end{equation}
The identities \eqref{d:pkn}, \eqref{idg:d}, \eqref{idg:B} and \eqref{idg:G}, 
Theorem~1 and a straightforward computation yield
\begin{itemize}
\item[1')]
$\SB (G_{2,4}, \gs_{2,4}) \in \{ 5, \dots, 9 \}$, \\
$\SB (G_{2,5}, \gs_{2,5}) \in \{ 7, \dots, 13 \}$, \\
$\SB (G_{2,6}, \gs_{2,6}) \in \{ 15, 16, 17\}$, \\
$\SB (G_{2,n}, \gs_{2,n}) = p_{2,n}+1$ for all $n \ge 7$,
\vspace{0.2em}
\item[2')]
$\SB (G_{k,n}, \gs_{k,n}) = p_{k,n}+1$ for all $k \ge 3$.
\end{itemize}
Corollary~\ref{c:grassmann} now follows together with the estimate
$
\SB (G_{k,n}, \gs_{k,n}) \, \le \, {n \choose k} ,
$
which is obtained by generalizing the embeddings~\eqref{e:chartCPn}
to $n \choose k$ symplectic embeddings $B^{2k(n-k)}(1) \to G_{k,n}$
covering $G_{k,n}$, see Lemma~4.1 and Section~6 in~\cite{Lu}.
\proofend

\end{document}